\documentclass[11pt]{amsart}

\usepackage{amscd, amsmath, amssymb, amsthm}
\usepackage[all]{xy}
\usepackage[dvips]{graphicx}
\usepackage[dvips]{color}

\theoremstyle{plain}
\newtheorem{Thm}{Theorem}[section]

\newtheorem{Cor}[Thm]{Corollary}

\newtheorem{Lem}[Thm]{Lemma}

\newtheorem{Prop}[Thm]{Proposition}

\newtheorem{Ques}[Thm]{Question}

\theoremstyle{remark}
\newtheorem{Exa}[Thm]{Example}
\newtheorem{Rem}[Thm]{Remark}

\theoremstyle{definition}
\newtheorem{Def}[Thm]{Definition}

\newcommand{\R}{\mathbb{R}}
\newcommand{\Z}{\mathbb{Z}}

\newcommand{\calA}{\mathcal{A}}

\newcommand{\calG}{\mathcal{G}}
\newcommand{\calI}{\mathcal{I}}

\newcommand{\abs}[1]{\lvert {#1} \rvert}
\newcommand{\Aut}{\mathrm{Aut}\,}

\newcommand{\emb}[2]{\mathrm{Emb}(\R^{#2} ,\R^{#1})}

\newcommand{\femb}[2]{\overline{\mathrm{Emb}} (\R^{#2} ,\R^{#1})}

\newcommand{\id}{\mathrm{id}}
\newcommand{\Imm}{\mathrm{Imm}}
\newcommand{\Int}{\mathrm{Int}\,}

\newcommand{\ord}{\mathrm{ord}\,}
\newcommand{\ori}{\mathrm{or}}

\newcommand{\pair}[2]{\langle #1 ,\, #2 \rangle}

\numberwithin{equation}{section}
\numberwithin{figure}{section}
\numberwithin{table}{section}

\begin{document}
\title{1-loop graphs and configuration space integral for embedding spaces}
\author[K. Sakai]{Keiichi Sakai}
\author[T. Watanabe]{Tadayuki Watanabe}
\address{Graduate School of Mathematical Sciences, The University of Tokyo, 3-8-1 Komaba, Meguro-ku, Tokyo 153-8914, Japan}
\address{Department of Mathematics, Hokkaido University, Kita 10, Nishi 8, Kita-Ku, Sapporo, Hokkaido, 060-0810, Japan}
\email{}
\thanks{KS is partially supported by (a) Grant-in-Aid for Young Scientists (B) 21740038, JSPS, (b) Grant for Basic Science Research Projects, the Sumitomo Foundation, and (c) The Iwanami Fujukai Foundation.
TW is partially supported by (a) Grant-in-Aid for JSPS Fellows 08J01880, JSPS, (b) Grants-in-Aid for Young Scientists (Start-up) 21840002, JSPS, (c) Start-up research funds, Hokkaido University.}

\date{\today}

\begin{abstract}
\def\baselinestretch{1.07}\footnotesize
We will construct differential forms on the embedding spaces $\emb{n}{j}$ for $n-j\ge 2$ using configuration space integral associated with 1-loop graphs, and show that some linear combinations of these forms are closed in some dimensions.
There are other dimensions in which we can show the closedness if we replace $\emb{n}{j}$ by $\femb{n}{j}$, the homotopy fiber of the inclusion $\emb{n}{j}\hookrightarrow\Imm (\R^j,\R^n)$.
We also show that the closed forms obtained give rise to nontrivial cohomology classes, evaluating them on some cycles of $\emb{n}{j}$ and $\femb{n}{j}$.
In particular we obtain nontrivial cohomology classes (for example, in $H^3(\emb{5}{2})$) of higher degrees than those of the first nonvanishing homotopy groups.
\end{abstract}
\maketitle


\setlength{\parskip}{2mm}

\section{Introduction}\label{sec_intro}
A {\em long immersion} is a smooth immersion $f: \R^j\to \R^n$ for some $n>j>0$ which agrees with the standard inclusion $\R^j\subset \R^n$ outside a disk $D^j\subset \R^j$.
A {\em long embedding} is an embedding $\R^j\hookrightarrow \R^n$ which is also a long immersion.
Let $\Imm(\R^j,\R^n)$ and $\emb{n}{j}$ be the spaces of long immersions and long embeddings respectively, both equipped with the $C^\infty$-topology.
In this paper we will construct some nontrivial cohomology classes of $\emb{n}{j}$ given by means of graphs.

Some graphs have appeared in previous works.
Some special graphs are introduced in \cite{R,CR} for describing a perturbative expansion of the BF theory functional integral for higher-dimensional embeddings, and an isotopy invariant of codimension two higher-dimensional embeddings is constructed via {\em configuration space integral} (CSI for short).
The graphs used in \cite{R,CR} are {\em 1-loop} graphs, i.e., those of the first Betti number exactly one (see also \cite{Wa}).

Recently Arone and Turchin announced that,
at least in the stable range $n\ge 2j+2$, the rational homology of $\emb{n}{j}$ can be expressed as the homology of some graph complex (see also \cite{ALV,To}).
On the other hand, a recent paper \cite{S} of the first author formally explains the invariance of the invariants of \cite{R,CR,Wa}
(in the cases when $n-j=2$) in the context of complexes of general graphs, which contain the graphs of \cite{R,CR}.
When the codimension $n-j$ is odd, a `$0$-loop'
graph cocycle of the complex of \cite{S}
gives the first nontrivial cohomology class of $\emb{n}{j}$ via CSI, which detects the lowest degree nontrivial homotopy class of $\emb{n}{j}$ given in \cite{B} (in odd codimension case).
These facts suggest that the method of graphs and CSI is effective even in the range $n<2j+2$.

In this paper we will focus on the 1-loop graphs of \cite{S} (which will be reviewed in \S\ref{sec:graph}).
We will construct some closed differential forms $z_k$ (resp.\ $\hat{z}_k$) of $\emb{n}{j}$ (resp.\ $\femb{n}{j}$) via CSI for arbitrary $n,j$ with $n-j\geq 2$ (see Theorems~\ref{thm:closed}, \ref{thm:closed2}).
Here $\femb{n}{j}$ is the homotopy fiber of $\emb{n}{j}\hookrightarrow\Imm (\R^j,\R^n)$ over the standard inclusion $\iota :\R^j\subset\R^n$.
Namely, $\femb{n}{j}$ is the space of smooth 1-parameter families of long immersions $\varphi_t:\R^j\to \R^n$, $t\in [0,1]$, such that
$\varphi_0=\iota$ and such that $\varphi_1\in \emb{n}{j}$.
The forgetting map 
\[ r: \femb{n}{j} \to \emb{n}{j} \]
given by $\{\varphi_t\}\mapsto \varphi_1$ is a fibration with homotopy fiber $\Omega\Imm(\R^j,\R^n)$.
The homotopy type of $\Imm(\R^j,\R^n)$ is well-known by \cite{Sm}.
So it follows that there is no big difference between the rational homotopy groups of $\emb{n}{j}$ and of $\femb{n}{j}$.

We will generalize the framework given in \cite{R, CR} to construct $z_k$ and $\hat{z}_k$.
They will be given explicitly as closed forms with values in $\calA_k=\calA_k(n,j)$, a vector space spanned by some graphs and quotiented by some diagrammatic relations (IHX/STU relations; see \S\ref{sec:graph}).
These forms represent nontrivial cohomology classes of $\emb{n}{j}$ and $\femb{n}{j}$ in dimensions stated in the following Theorem.

\begin{Thm}[Theorems~\ref{thm:closed}, \ref{thm:closed2}, \ref{thm:nontrivial}]
$H^{(n-j-2)k}_{DR}(\emb{n}{j};\calA_k)$ is nontrivial if one of the following holds and if $k\ge 2$ is such that the space $\calA_k=\calA_k(n,j)$ does not vanish (see Proposition \ref{prop:calA} below):
\begin{itemize}
\item $n$ is odd.
\item $n$ is even, $j$ is odd, and $k\le 4$.
\item $n\ge 12$ is even and $j=3$.
\item $n,j$ are both even, $n-j>2$ and $k$ is large enough so that $2k(n-j-2)>j(2n-3j-3)$.
\end{itemize}
$H^{(n-j-2)k}_{DR}(\femb{n}{j};\calA_k)$ is nontrivial if both $n,j$ are even and if $k$ is such that $\calA_k\ne 0$.
See Figure~\ref{fig:range_j_n}.
\end{Thm}

\begin{Prop}[\S\ref{subsec:n-j=even}, Proposition~\ref{prop:A3_odd_codim}]\label{prop:calA}
In even codimension case, $\calA_k\cong\R$ if $k\not\equiv n$ modulo $2$, and $\calA_k=0$ otherwise.
When $n$ is odd and $j$ is even, $\calA_3\cong\R$.
\end{Prop}

When one of $n$ and $j$ is odd, the cohomology class $[z_k]$
generalizes
invariants of \cite{R,CR,Wa} for codimension two long embeddings in $\R^n$, which can be regarded as in $H^0_{DR}(\emb{n}{n-2})$.
All of our cohomology classes are of higher degrees than those discussed in \cite{B} and hence new.

The construction of the closed forms $z_k$ and $\hat{z}_k$ will be given in \S\ref{s:closed-form}.
For this, we need the following extra arguments in addition to those of \cite{R, CR}.
\begin{enumerate}
\item
 In even codimension case, we need lemmas of \cite{S} (in addition to those of \cite{R, CR}) to show the vanishing of the obstructions to the closedness which arise from degenerations of certain kind of subgraphs.
\item
 In odd codimension case, we should take more general 1-loop graphs \cite{S} than those in \cite{R, CR} into account in order to get meaningful closed forms.
Moreover, we will generalize the cancellation arguments due to the diagrammatic relations to those of more general kinds of subgraph degenerations.
\item
 In the case when both $n,j$ are even, almost all the obstructions as above cancel, but we have no proof of the vanishing of so-called `anomaly' arising from degenerations of whole graphs.
So we consider another space $\femb{n}{j}$ on which we can construct a correction term.
See \S\ref{ss:anomaly_correction}.

In fact the correction term restricts to a closed form of $\Omega\Imm(\R^j,\R^n)$.
It seems likely that this closed form is related to the surjection $\pi_0(\emb{5}{3})\to 24\Z$ given by Smale-Hirsch map \cite{Ek,HM}.
See Remark~\ref{rem:anomaly-SH}.
\end{enumerate}


To prove the nontriviality of $[z_k]$ and $[\hat{z}_k]$, we will generalize in \S\ref{sec_ribbon} the method of \cite{Wa} to higher-dimensions to construct nontrivial homology classes of $\emb{n}{j}$ and $\femb{n}{j}$ by a `resolution of crossings', an analogous technique to that considered in \cite{CCL}.
We will explicitly coompute the pairings of these homology classes with $z_k$ and $\hat{z}_k$, and show that they are not zero.

There seems to be further possible progress in the direction of this paper.
The nontriviality results of this paper might be generalized for graphs with one or more loop components, if the corresponding forms were proved to be closed.
There might be other generalizations as in \cite{Wa2}.
Indeed, some cycles on $\emb{2k+1}{k}$ are constructed in \cite{Wa2}, which can be considered as a generalization of the construction of this paper.
It would be also interesting to ask how our cohomology classes given in terms of graphs relate to the actions of little cubes operad \cite{B1}.

Now we give an account how the authors began writing this paper.
A part of the present paper is based on a note by KS (the first author).
After \cite{Wa} has been published KS arrived at the result of the present paper for both $n$ and $j$ odd, and wrote the detailed proof into a note.
But TW (the second author) had a proof of the same result independently and in fact, after his note has been written KS was informed about the preprint of TW in which a rough sketch of the same result is given.
So the authors decided to work together and extended the main result of the note to arbitrary pairs $n>j\geq 2$, $n-j\geq 2$.

\section*{Acknowledgment}
The authors express their great appreciation to Professor Toshitake Kohno for his encouragement to the authors.
The authors are also grateful to Ryan Budney for informing the authors about his result on the connectivity of $\emb{n}{j}$, to Masamichi Takase for useful discussions and suggestions, and to the referee for suggesting possible improvements of the earlier version of this paper.

\section{1-loop graphs}\label{sec:graph}

In this section we review the definition of graphs introduced in \cite{S}, which generalize those appearing in \cite{CR, R, Wa}.

\subsection{Graphs}\label{subsec:graphs}

A {\em graph} in this paper has two kinds of vertices, namely {\em external} vertices (or shortly {\em e-vertices}) and {\em internal} ones (shortly {\em i-vertices}), and two kinds of edges, $\theta${\em -edges} and $\eta${\em -edges}.
We depict e- and i-vertices as $\circ$ and $\bullet$ respectively.
We depict $\theta$-edges and $\eta$-edges as dotted lines and solid lines, respectively.
We assume that no single edge forms a loop.

\begin{Def}\label{def:admissible}
A vertex $v$ of a graph is said to be {\em admissible} if it is at most trivalent and is one of the following forms;
\[
 \includegraphics{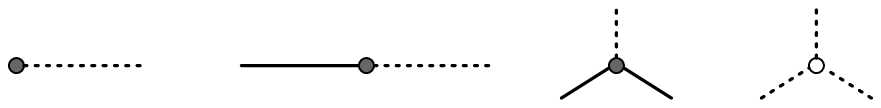}
\]
A graph is said to be {\em admissible} if all its vertices are admissible.\qed
\end{Def}

\begin{Rem}
By definition, the endpoints of an $\eta$-edge of an admissible graph must be i-vertices.
Those of a $\theta$-edge can be either i- or e-vertices.
In \cite{S} the vertices shown in Definition \ref{def:admissible} were said to be admissible and `non-degenerate'.\qed
\end{Rem}

\begin{Def}\label{def:order}
Below {\em 1-loop graph} means an admissible graph whose first Betti number is one.
The {\em order} of a 1-loop graph $\Gamma$, denoted by $\ord (\Gamma )$, is half the number of the vertices of
$\Gamma$ ($\ord (\Gamma )$ is a positive integer; see Remark \ref{rem:N}).\qed
\end{Def}

\begin{Exa}\label{ex:1-loop_graph}
The following three graphs are examples of admissible 1-loop graphs.
\[
 \includegraphics{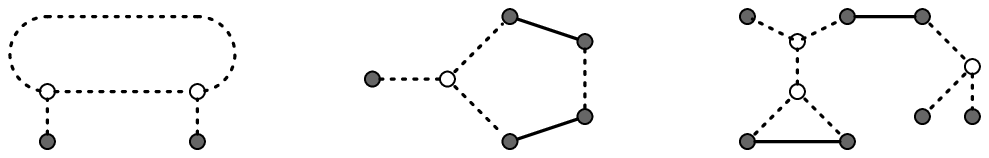}
\]
The orders of these graphs are $2$, $3$ and $5$ respectively.
A graph may have a large tree subgraph which shares only one vertex with the unique cycle, like the third graph
(such graphs have not been considered in \cite{CR,R}).\qed
\end{Exa}

\subsection{Labels and orientations of graphs}\label{subsec:ori_graph}

Below let $(n,j)$ be a pair of positive integers with $n-j \ge 2$.
Here we introduce the notion of {\em labelled graphs}.

\begin{Def}\label{def:labelled_graph}
Denote by $V_i(\Gamma )$, $V_e (\Gamma )$, $E_{\eta}(\Gamma )$ and $E_{\theta}(\Gamma )$ the sets of all i-vertices,
e-vertices, $\eta$-edges and $\theta$-edges of a graph $\Gamma$, respectively.
We also write $V(\Gamma ):= V_i (\Gamma )\cup V_e (\Gamma )$ and
$E(\Gamma ):= E_{\eta}(\Gamma )\cup E_{\theta}(\Gamma )$.
We decompose $V(\Gamma ) \sqcup E(\Gamma )$ into two disjoint subsets $S(\Gamma )$ and $T(\Gamma )$ given by
\[
 (S(\Gamma ) , T(\Gamma )) :=
 \begin{cases}
  (V(\Gamma ), E(\Gamma )) & n,j \text{ odd}, \\
  (E(\Gamma ), V(\Gamma )) & n,j \text{ even}, \\
  (V_e (\Gamma )\cup E_{\eta}(\Gamma ),V_i(\Gamma )\cup E_{\theta}(\Gamma )) & n \text{ odd, } j \text{ even}, \\
  (V_i (\Gamma )\cup E_{\theta}(\Gamma ),V_e(\Gamma )\cup E_{\eta}(\Gamma )) & n \text{ even, } j \text{ odd}.
 \end{cases}
\]
Below we will write $k_S := \abs{S(\Gamma )}$ and $k_T := \abs{T(\Gamma )}$.
A {\em labelled graph} is a 1-loop, admissible graph $\Gamma$ together with bijections
\[
 \rho_1:\{ 1,\dots ,k_S\}\longrightarrow S(\Gamma ),\quad
 \rho_0:\{ 1,\dots ,k_T\}\longrightarrow T(\Gamma ).\qed
\]
\end{Def}

\begin{Rem}\label{rem:N}
It holds $2\abs{E_{\theta}(\Gamma )}-3\abs{V_e (\Gamma )}-\abs{V_i (\Gamma )}=0$ since exactly one (resp.\ three)
$\theta$-edge(s) emanates from each i-vertex (resp.\ e-vertex).
Hence $\abs{V_e (\Gamma )}+\abs{V_i (\Gamma )}=2\abs{E_{\theta}(\Gamma )}-2\abs{V_e (\Gamma )}$. 
This implies that $\ord (\Gamma )$ is an integer
and is equal to $\abs{E_{\theta}(\Gamma )}-\abs{V_e (\Gamma )}$
(in \cite{S} the order was defined as the latter number).
Putting $k:=\ord (\Gamma )$, we can show that $k_S =k_T =2k$ in even codimension case, and
$(k_S ,k_T )=(3k,k)$ ($n$ odd, $j$ even) or $(k,3k)$ ($n$ even, $j$ odd).\qed
\end{Rem}

To fix the signs of the configuration space integrals (see \S\ref{s:closed-form}), we orient the graphs following
\cite[Appendix B]{T} so that the elements of $S(\Gamma )$ (resp.\ $T(\Gamma )$) are of odd (resp.\ even) degrees.

\begin{Def}\label{def:ori_graph}
We think of an edge $e$ as a union of two shorter segments; $e=h_1(e)\cup h_2(e)$, $h_1(e)\cap h_2(e)=\text{the midpoint of }e$.
Each $h_i(e)$ is called a {\em half-edge} of $e$.

For an edge $e$, define $H(e)=\{ h_1(e),h_2(e)\}$ as the set of half-edges of $e$.
For any graph $\Gamma$, define a graded vector space $Ori(\Gamma )$ by
\[
 Ori (\Gamma ) := \R S(\Gamma ) \oplus \R T(\Gamma ) \oplus \bigoplus_{e\in E(\Gamma )}\R H(e),
\]
here 
$\R X:=\bigoplus_{x\in X}\R x$ for a set $X$, and we regard $Ori(\Gamma )$ as a graded vector space by assigning the degrees to the elements of $S(\Gamma )$, $T(\Gamma )$ and $H(e)$ as in Table~\ref{tab:degrees}.
\begin{table}[htb]
\begin{center}
\begin{tabular}{|c|c|c|c|c|c|}
\hline
i-vertices & e-vertices & $\eta$-edges & $\theta$-edges & half $\eta$-edges & half $\theta$-edges \\
\hline
$j$        & $n$        & $j-1$        & $n-1$          & $j$               & $n$ \\
\hline
\end{tabular}
\end{center}
\caption{Degrees of elements of $Ori(\Gamma )$}
\label{tab:degrees}
\end{table}
An {\em orientation} of a graph $\Gamma$ is that of one dimensional vector space $\det Ori(\Gamma )$, where $\det V :=\bigwedge^{\dim V}V$ for a vector space $V$.
An orientation of a labelled graph $\Gamma$ is determined by its edge-orientation.
We denote an orientation determined in this way by $o =\ori (\Gamma )$.\qed
\end{Def}

See \S\ref{subsec:integral} for the meaning of the labelled graphs and their orientations as above.

\subsection{A graph cocycle}\label{subsec:graph_cocycles}

\begin{figure}[htb]
\includegraphics{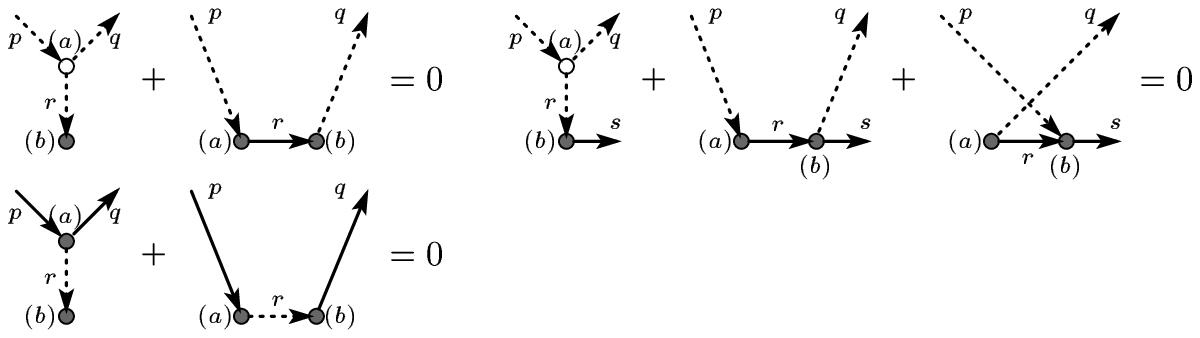}
\caption{ST, ST2 and C relations, even codimension case}\label{fig:relations_even}
\end{figure}

\begin{figure}[htb]
\includegraphics{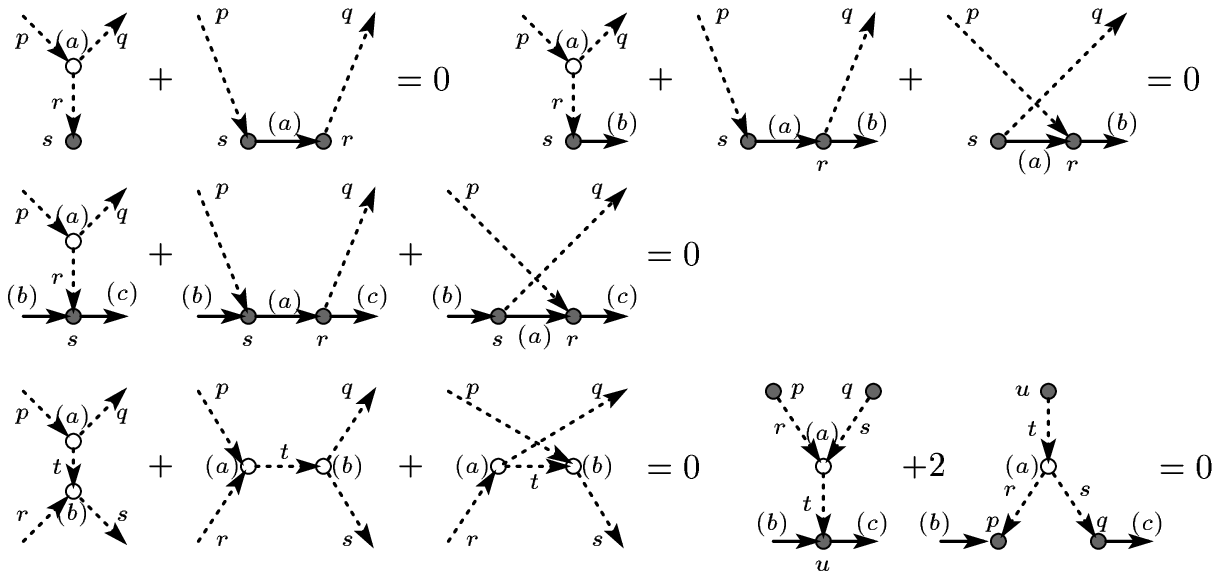}
\caption{ST, ST2, STU, IHX and Y relations, odd codimension case}\label{fig:relations_odd}
\end{figure}

\begin{figure}[htb]
\includegraphics{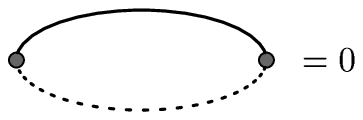}
\caption{L relation (for arbitrary $n$ and $j$)}\label{fig:relation_l}
\end{figure}

\begin{Def}\label{def:relation_graph}
Denote by $\tilde{\calG}_k =\tilde{\calG}_k (n,j)$ the set of labelled, oriented 1-loop graphs $(\Gamma ,\ori (\Gamma ))$
of order $k$ (the definitions of labels and orientations depend on the parities of $n,j$).
Define the vector space $\calG_k =\calG_k (n,j)$ of labelled, oriented graphs by
\[
 \calG_k := \R \tilde{\calG}_k /(\Gamma , -\ori (\Gamma )) \sim -(\Gamma ,\ori (\Gamma )),
\]
where $-\ori (\Gamma )$ is the orientation obtained by reversing the edge-orientation (that is, $\R H(e)$-part) of $\ori (\Gamma )$.
Define the vector space $\calA_k=\calA_k (n,j)$ by
\[
 \calA_k := \calG_k / \text{relations, labels}
\]
where relations are shown in Figures \ref{fig:relations_even}, \ref{fig:relations_odd} and \ref{fig:relation_l} and
the quotient by ``labels'' means that we regard two labelled oriented graphs with the same underlying oriented graphs
as being equal to each other in $\calA_k$.
Each $[\Gamma ] \in \calA_k$ possesses an orientation induced from $\ori (\Gamma )$ of $\Gamma \in \calG_k$.
In Figures \ref{fig:relations_even}, \ref{fig:relations_odd} and \ref{fig:relation_l}, we have already forgotten the
labels.
The orientations of graphs are indicated by the letters assigned to vertices and edges (which correspond to $\R S(\Gamma )\oplus\R T(\Gamma )$-part of $\ori (\Gamma)$), and the orientations of edges (which correspond to $\R H(e)$-part).
When $(a),(b),\dots$ are numbers for $S(\Gamma )$ (resp.\ $T(\Gamma )$), then $p,q,\dots$ are those for $T(\Gamma )$
(resp.\ $S(\Gamma )$).\qed
\end{Def}


\begin{Rem}
In \cite{S} we introduced `graph complexes,' whose coboundary operation $\delta$ is given as a signed sum of graphs obtained by contracting the edges one at a time (we have several complexes depending on the parities of $n$ and $j$).
We defined the relations in Figures \ref{fig:relations_even}, \ref{fig:relations_odd} and \ref{fig:relation_l} so that the linear combination
\begin{equation}\label{eq:X_k}
 X_k := \frac{1}{k_S ! k_T !} \sum_{\Gamma} [\Gamma ] \otimes \Gamma \in \calA_k \otimes \calG_k
\end{equation}
of graphs with (untwisted) coefficients in $\calA_k$, where the sum runs over all the labelled graphs of order $k$ (with an orientation
assigned), becomes a `cocycle', i.e., $\delta X_k =0$.
This vanishing is an algebraic expression of the cancellation of fiber integrations along the `principal faces' of the boundary of compactified configuration spaces; see \S\ref{s:outline-proof}. 

The Y relation is needed to construct cocycles in odd codimension case.
In $\calA_3$, the Y relation is a consequence of the STU and the IHX relations (but it might not hold for general
$\calA_k$).\qed
\end{Rem}

In \S\ref{s:closed-form} closed forms of $\emb{n}{j}$ (or $\femb{n}{j}$) with coefficients in $\calA_k$ will be defined.
What we know about $\calA_k$ are stated in Proposition~\ref{prop:calA} and will be proved in \S\ref{sec:calA}.



\section{Cohomology classes of embedding spaces from configuration space integral}\label{s:closed-form}

\newcommand{\fig}[1]
        {\raisebox{-0.5\height}
                 {\includegraphics{#1}}
        }
\newcommand{\tfig}[1]
		{\raisebox{-0.4\height}
				{\input{#1}}
		}
\def\red#1{{\color{red}#1}}
\def\blue#1{{\color{blue}#1}}
\setlength{\parskip}{2mm}

\subsection{Configuration space integral}\label{subsec:integral}

Let $\varphi:\R^j\hookrightarrow \R^n$ denote a long embedding. Let $\Gamma=(\Gamma ,\ori )$ be an oriented graph with $s$ i-vertices and $t$ e-vertices\red{,} labelled by the bijections $\rho_1$ and $\rho_0$ (\S\ref{subsec:ori_graph}). Then consider the space
\[
 C^o_{\Gamma} := \{
 (\varphi;x_1 ,\dots ,x_s ; x_{s+1} ,\dots ,x_{s+t}) \in \emb{n}{j} \times C^o_s (\R^j ) \times C^o_t (\R^n ) \, | \,
 \varphi(x_p ) \ne x_{s+q} ,\ \forall p, q>0
 \},
\]
where $C^o_k (M)$ denotes the configuration space in the usual sense;
\[
 C_k^o (M) := \{ (x_1 ,\dots ,x_k ) \in M^{\times k} \, | \, x_i \ne x_j \text{ if } i \ne j \}.
\]
The space $C^o_\Gamma$ is naturally fibered over $\emb{n}{j}$, namely, the projection map 
\[
 \pi_{\Gamma} : C^o_{\Gamma} \longrightarrow \emb{n}{j},
\]
given by $(\varphi;x_1 ,\dots ,x_s ; x_{s+1} ,\dots ,x_{s+t})\mapsto \varphi$, is a fiber bundle with fiber $C^o_\Gamma(\varphi)=C^o_s (\R^j ) \times C^o_t (\R^n )\setminus \bigcup_{{1\leq i\leq s}\atop{s+1\leq j\leq s+t}}\{\varphi(x_i)=x_j\}$.

From now on we will define for each oriented graph $\Gamma$ a differential form $I(\Gamma)$ on $\emb{n}{j}$ as the fiber integral of the following form
\[ I(\Gamma)=\pm (\pi_\Gamma)_*\bigwedge_{e\in E(\Gamma)} \omega_e. \]
Here $(\pi_\Gamma)_*:\Omega_{DR}^{*}(C_\Gamma^o)\to\Omega_{DR}^*(\emb{n}{j})$ denotes the integration along the fiber, $\omega_e$ is the `edge form' (see below for precise definition). The choice of a sign from a graph orientation will make the definition rather complicated.

Precise definition of $I(\Gamma)$ is as follows. The bijections $\rho_{1}$ and $\rho_0$ give an orientation 
\[ \mathrm{or}'(\Gamma):=\rho_{1}(1)\wedge\cdots\wedge\rho_{1}(k_S)\wedge\rho_{0}(1)\wedge\cdots\wedge\rho_0(k_T)\]
 of $\R{S(\Gamma)}\oplus \R{T(\Gamma)}$. We arrange $\mathrm{or}'(\Gamma)$ in the form $(\mbox{i-vertices})\wedge(\mbox{e-vertices})\wedge(\mbox{$\eta$-edges})\wedge(\mbox{$\theta$-edges})$ as 
\begin{equation}\label{eq:or_order}
\begin{split}
 \mathrm{or}'(\Gamma)=\varepsilon(\rho_1,\rho_0)
	\bigwedge_{p=1}^s \rho_{\underline{j}}(i_p)\wedge
	\bigwedge_{q=1}^t \rho_{\underline{n}}(j_q)\wedge
	\bigwedge_{r=1}^{|E_\eta(\Gamma)|} \rho_{\underline{j-1}}(\sigma_r)\wedge
	\bigwedge_{u=1}^{|E_\theta(\Gamma)|} \rho_{\underline{n-1}}(\tau_u) 
\end{split}
\end{equation}
for $\varepsilon(\rho_1,\rho_0)=\pm 1$, $i_1<\cdots<i_s$, $j_1<\cdots<j_t$ and for some numbers $\sigma_r$, $\tau_u$, which are uniquely chosen up to even swappings. Here $\underline{p}$ denotes $p\ \mathrm{mod}\ 2$. The vertex part of (\ref{eq:or_order}) determines a bijection 
\[ v:  V(\Gamma)\to \{1,\ldots,s+t\}\quad{\mbox{by}}\]
\[ v^{-1}(p)=\left\{
	\begin{array}{ll}
		\rho_{\underline{j}}(i_p) & \mbox{if $1\leq p\leq s$}\\
		\rho_{\underline{n}}(j_{p-s}) & \mbox{if $s+1\leq p\leq s+t$}
	\end{array}\right. \]
Now we orient edges of $\Gamma$ so that $\mathrm{or}'(\Gamma)$ and the edge orientation give the orientation $\mathrm{or}(\Gamma)$ where an arrow $\overrightarrow{ab}$ on an edge $ab$ from a vertex $a$ to a vertex $b$ corresponds to $h_a\wedge h_b\in\mathrm{\det}\R H(ab)$ of the half edges $h_a, h_b$ including $a, b$ respectively. To each oriented edge $e=\overrightarrow{ab}$ of $\Gamma$, we assign a map
$ \phi_e : C^o_{\Gamma} \longrightarrow S^{N-1}$ where $N=j$ or $n$ according to whether $e$ is an $\eta$- or a $\theta$-edge, defined by
\[
 \phi_e (\varphi;x_1,\ldots,x_s;x_{s+1},\ldots,x_{s+t}) :=
		\frac{z_{v(b)} -z_{v(a)}}{\abs{z_{v(b)} -z_{v(a)}}}
\quad\mbox{where}\]
\[
 z_{v(p)} :=
 \begin{cases}
  x_{v(p)}     & \text{if } e \text{ is an } \eta \text{-edge (and hence } a,b \text{ are both i-vertices)}, \\
       & \text{\ \  or if } e \text{ is a } \theta \text{-edge and } p \text{ is an e-vertex},\\
  \varphi(x_{v(p)} ) & \text{if } e \text{ is a } \theta \text{-edge and } p \text{ is an i-vertex}. \\
 \end{cases}
\]
Let $vol_{S^{N-1}}$ denote the volume form of $S^{N-1}$ which is (anti)symmetric with respect to the antipodal map $\Upsilon:S^{N-1}\to S^{N-1}$, i.e. $\Upsilon^* vol_{S^{N-1}} = (-1)^N vol_{S^{N-1}}$, and is normalized as $\int_{S^{N-1}}vol_{S^{N-1}}=1$, and define the `edge form' by
\[
 \omega_e :=\phi^*_e vol_{S^{N-1}} \in \Omega^{N-1}_{DR}(C^o_{\Gamma}).
\]

We define $\omega_\Gamma\in\Omega^*_{DR}(C^o_{\Gamma})$ by
\begin{equation}\label{eq:omega_G}
 \omega_{\Gamma} := \varepsilon(\rho_1,\rho_0)
	\bigwedge_{r=1}^{|E_\eta(\Gamma)|}\omega_{\rho_{\underline{j-1}}(\sigma_r)}
	\wedge\bigwedge_{u=1}^{|E_\theta(\Gamma)|}\omega_{\rho_{\underline{n-1}}(\tau_u)}.
\end{equation}
The integration of $\omega_{\Gamma}$ along the fiber of the bundle $\pi_{\Gamma}$ given above yields a differential form on $\emb{n}{j}$;
\[
 I(\Gamma ) := (\pi_{\Gamma})_* \omega_{\Gamma} \in \Omega^{*}_{DR}(\emb{n}{j}).
\]
Here the orientation on the fiber is imposed by the canonical one given by $dx_1\wedge\cdots\wedge dx_{s+t}$, $dx_i=dx_i^{(1)}\wedge\cdots\wedge dx_i^{(N)}$, $N=n$ or $j$. If $\Gamma$ is an admissible 1-loop graph of order $k$, then the degree of $I(\Gamma )$ is
$(n-j-2)k$ (see \cite{S}).

\begin{Prop}
The integral $I(\Gamma)$ converges. So we have a well-defined linear map $I:\calG_k \to \Omega^{(n-j-2)k}_{DR}(\emb{n}{j})$.\qed
\end{Prop}

\begin{Rem}\label{rem:boundary-conf}
Since the fiber $C^o_\Gamma(\varphi)$ of $\pi_{\Gamma}$ is not compact, the convergence of the integral is not trivial. As was done in \cite{BT, R}, the proof of the convergence uses a compactification $C_\Gamma(\varphi)$ of $C_\Gamma^o(\varphi)$, obtained by `blowing-up' along the stratification formed by all the singular strata in the product $\varphi(S^j)^{\times s}\times (S^{n})^{\times t}$ where some points come close to each other or go to infinity. Here we identify $\R^j$ (resp. $\R^n$) with the complement of a point $\infty$ in $S^j$ (resp. $S^n$) and $\varphi$ extends uniquely and smoothly to $S^j$ by mapping $\infty$ to $\infty$. The result of the blow-ups is a smooth manifold with corners, stratified by possible parenthesizations of $s+t$ distinct letters corresponding to the $s+t$ points. The parenthesis corresponds to a degeneration of the parenthesized points collapsed into a multiple point. In particular, the codimension one (boundary) strata is given by a word with one pair of parentheses which encloses a subset $A\subset V(\Gamma)\cup\{\infty\}$. Note that the resulting manifold with corners depends only on $\varphi$ and the numbers $(s,t)$. In the case where $s=0$, we will denote the result by $C_t(\R^n)$ and in the case where $t=0$, we will denote the result by $C_s(\R^j)$. See for example \cite{BT, R} for detail of the compactification.\qed
\end{Rem}

Now we define the main differential form of this paper:
\[ z_k:=(1\otimes I)(X_k)\in \calA_k\otimes \Omega^{(n-j-2)k}_{DR}(\emb{n}{j}) \]
where $X_k\in \calA_k\otimes \calG_k$ is defined in (\ref{eq:X_k}). 

We will see that the differential form $z_k$ is closed for approximately half of the pairs $(n,j)$ with $n-j\geq 2$. However we do not know whether $z_k$ is closed for all $(n,j)$ due to some `anomaly'. When the anomaly may exist we consider the pullback of $z_k$ to $\femb{n}{j}$ and we will introduce (in \S\ref{ss:anomaly_correction}) a correction term $\Theta_k\in \calA_k\otimes\Omega_{DR}^{(n-j-2)k}(\femb{n}{j})$ for the anomaly and define
\begin{equation}\label{eq:hatz}
 \hat{z}_k:=r^*z_k-\Theta_k\in \calA_k\otimes\Omega_{DR}^{(n-j-2)k}(\femb{n}{j}). 
\end{equation}

\begin{figure}
\fig{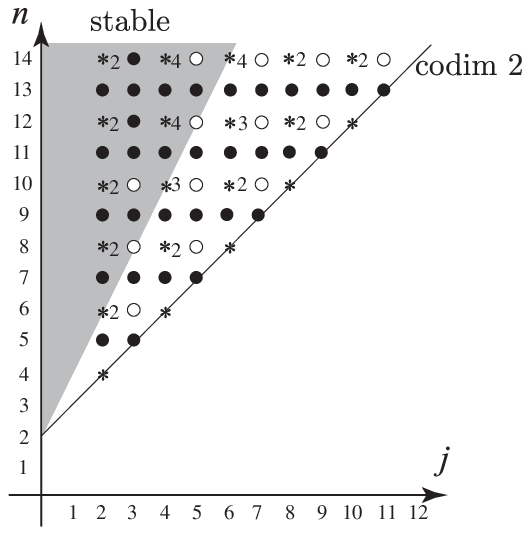}
\caption{$\bullet$ is a pair of dimension $(j,n)$ where $z_k\in\calA_k\otimes\Omega_{DR}^{(n-j-2)k}(\emb{n}{j})$ is proved to be closed for all $k\geq 2$, $\circ$ is a pair $(j,n)$ where $z_k$ is proved to be closed for $k\leq 4$, $*$ is a pair where $\hat{z}_k\in\calA_k\otimes\Omega_{DR}^{(n-j-2)k}(\femb{n}{j})$ is proved to be closed for all $k\geq 2$, and $* p$ indicates that $\hat{z}_k$ descends to the closed form $\bar{z}_k$ on $\emb{n}{j}$ for all $k\geq p$. We will show in \S\ref{sec_ribbon} that $z_k$ or $\hat{z}_k$ in the range shown in this figure are nontrivial, provided that $\calA_k\neq 0$.}
\label{fig:range_j_n}
\end{figure}

\begin{Thm}\label{thm:closed}
Let $n, j, k$ be positive integers with $n-j\geq 2$, $n\geq 4$, $k\geq 2$.
\begin{enumerate}
\item The form $z_k\in \calA_k\otimes\Omega_{DR}^{(n-j-2)k}(\emb{n}{j})$ is closed if one of the following holds:
\begin{enumerate}
\item $n$: odd ($j$ may be both odd and even).
\item $n$: even, $j$: odd, $k\leq 4$.
\item $n$: even $\geq 12$, $j=3$.
\end{enumerate}
(see Figure~\ref{fig:range_j_n}, $\bullet$ and $\circ$).
\item The form $\hat{z}_k\in \calA_k\otimes\Omega_{DR}^{(n-j-2)k}(\femb{n}{j})$ is closed if both $n$ and $j$ are even. (See Figure~\ref{fig:range_j_n},  $*$).
\end{enumerate}
\end{Thm}

Theorem~\ref{thm:closed} generalizes a result of \cite{CR}, which is concerned with the cases (1) $n$, $j$: odd, $n=j+2$, (2) $(n,j,k)=(4,2,3)$. The correction term for the latter case considered in \cite{CR} is different from ours but their invariant is well-defined on $\emb{4}{2}$.

\begin{Thm}\label{thm:closed2}
If $n-j>2$, $n,j$ both even and $k > \dfrac{j(2n-3j-3)}{2(n-j-2)}$, then there exists an $((n-j-2)k+j)$-form $\bar\alpha_k$ on $C_1(\R^j)\times \emb{n}{j}$ such that the form 
\[ \bar{z}_k:=z_k-\int_{C_1(\R^j)}\bar\alpha_k\in \calA_k\otimes\Omega_{DR}^{(n-j-2)k}(\emb{n}{j})\]
where $\int_{C_1(\R^j)}$ denotes the integration along the fiber, is closed and that its pullback to $\femb{n}{j}$ represents the same cohomology class as $\hat{z}_k$. (See Figure~\ref{fig:range_j_n},  $* p$).
\end{Thm}

\subsection{Outline}\label{s:outline-proof}
As usual in the theory of configuration space integral, the proof of Theorem~\ref{thm:closed} is reduced to the vanishing of integrals over the boundary of the fiber by the generalized Stokes theorem. Now we shall give a quick review of the necessary arguments in the proof, following \cite{R}. Recall that the generalized Stokes theorem for a fiber bundle $\pi:E\to B$ and a differential form $\alpha\in\Omega_{DR}^*(E)$ states that:
\begin{equation}\label{eq:stokes}
 d\pi_*\alpha=\pi_*d\alpha+J\pi_*^\partial \alpha,  
\end{equation}
where $J\gamma=(-1)^{\deg{\gamma}}\gamma$ and $\pi^\partial$ is $\pi$ restricted to the boundary of the fiber. Here the orientation of the boundary of the fiber is imposed by the inward-normal-first convention. Applying the generalized Stokes theorem (\ref{eq:stokes}) to $\pi_\Gamma$ we have
\begin{equation}\label{eq:dz}
 dz_k=\frac{1}{k_S!k_T!}\sum_{{\Gamma}\atop{\mathrm{labelled}}}[\Gamma]\otimes J(\pi_\Gamma)^\partial_*\omega_\Gamma
	=\frac{1}{k_S!k_T!}\sum_{{\Gamma}\atop{\mathrm{labelled}}}[\Gamma]\otimes J\sum_{A\subset V(\Gamma)}(\pi_\Gamma^{\partial_A})_*\omega_\Gamma. 
\end{equation}
Here $\pi^{\partial_A}_{\Gamma}$ is $\pi_{\Gamma}$ restricted to the codimension one face $\Sigma_A(\varphi)$ of $\partial C_{\Gamma}(\varphi)$ corresponding to the collapse of points in $A\subset V(\Gamma)$ (see Remark~\ref{rem:boundary-conf}). 

Each codimension one stratum $\Sigma_A$ is the pullback in the following commutative square:
\begin{equation}\label{eq:Sigma_A} \xymatrix{
	\Sigma_A \ar[r]^-{\hat{D}_A} \ar[d]_-{p_A} & \hat{B}_A \ar[d]^{\rho_A}\\
	C_{\Gamma/\Gamma_A} \ar[r]^-{D_A} & \calI_{j}(\R^n)
}\end{equation}
Here $\Gamma_A\subset \Gamma$ is the maximal subgraph with $V(\Gamma_A)=A$, $\Gamma / \Gamma_A$ is $\Gamma$ with the subgraph $\Gamma_A$ collapsed into a point.  Each term in the left hand vertical column of the square diagram is fibered $\Sigma_A=\Sigma_A(\emb{n}{j})$, $C_{\Gamma/\Gamma_A}=C_{\Gamma/\Gamma_A}(\emb{n}{j})$ over $\emb{n}{j}$, or over $\femb{n}{j}$ by the pullback along $r$. The right hand vertical column $\rho_A$ itself is a fiber bundle over $\calI_j(\R^n)$. The entries of the right hand vertical column of the diagram are given as follows: $\calI_j (\R^n )$ is the space of linear injective maps $\R^j \hookrightarrow \R^n$, the fiber $\hat{B}_A(f)$ of $\rho_A$ over $f\in \calI_j(\R^n)$ is the `microscopic' configuration space, i.e., $C_{\Gamma_A}(f)$ quotiented by the actions of overall translations of points along $f(\R^j)$ and overall dilations in $\R^n$ around the origin. Then the integral for $\Gamma$ restricted to the codimension one face $\Sigma_A$ is written as
\[
 (\pi^{\partial_A}_\Gamma)_*\omega_\Gamma
	=\int_{C_{\Gamma/\Gamma_A}}
		D_A^*\rho_{A*}\hat{\omega}_{\Gamma_A}\wedge \omega_{\Gamma/\Gamma_A}
\]
where $\int_{C_{\Gamma/\Gamma_A}}$ denotes the integration along the fiber, $\hat\omega_{\Gamma_A}\in\Omega_{DR}^*(\hat{B}_A)$ is the wedge of $\omega_e$'s for $\Gamma_A$ defined as in (\ref{eq:omega_G}). Note that $\deg\hat\omega_{\Gamma_A}=|E_\theta(\Gamma_A)|(n-1)+|E_\eta(\Gamma_A)|(j-1)$, $\deg \rho_{A*}\hat\omega_{\Gamma_A}=\deg\hat\omega_{\Gamma_A}-|V_e(\Gamma_A)|n-(|V_i(\Gamma_A)|-1)j+1$. 

With these facts in mind, the proof of Theorem~\ref{thm:closed} can be outlined as follows, which looks quite similar to that of the invariance proof of the invariant of \cite{R, CR} (but the detail is somewhat different).
\begin{proof}[Outline of the proof of Theorem~\ref{thm:closed}]
As in \cite{R}, the codimension one faces are classified into the following types, depending on the method of proof of vanishing of the integrals of (\ref{eq:dz}).
\begin{enumerate}
\item (Principal face) $\Sigma_A$ for $|A|=2$.
\item (Hidden face) $\Sigma_A$ for $2<|A|<|V(\Gamma)|$ corresponding to non-infinite diagonals.
\item (Infinite face) $\Sigma_A$ for $1\leq |A|\leq|V(\Gamma)|$ corresponding to diagonals involving the infinity.
\item (Anomalous face) $\Sigma_A$ for $A=V(\Gamma)$.
\end{enumerate}
In the sum (\ref{eq:dz}) the vanishing of the contribution of the principal faces has essentially been given a proof in \cite{S} in a general terms of the graph complex. But we give another explanation for the special cycle $X_k$ of the graph complex, namely, explain how the relations in \S\ref{sec:graph} work to prove the vanishing of the principal faces contributions. We only give here a proof of the vanishing given by the STU relation when $n$ is odd and $j$ is even because the other relations work similarly. 

Let $\Gamma_1,\ldots,\Gamma_6$ be as in Figure~\ref{fig:stu_labeled}. ($\Gamma_5$, $\Gamma_6$ are unnecessary if the bottom i-vertex of $\Gamma_1$ is univalent.)
\begin{figure}
\includegraphics{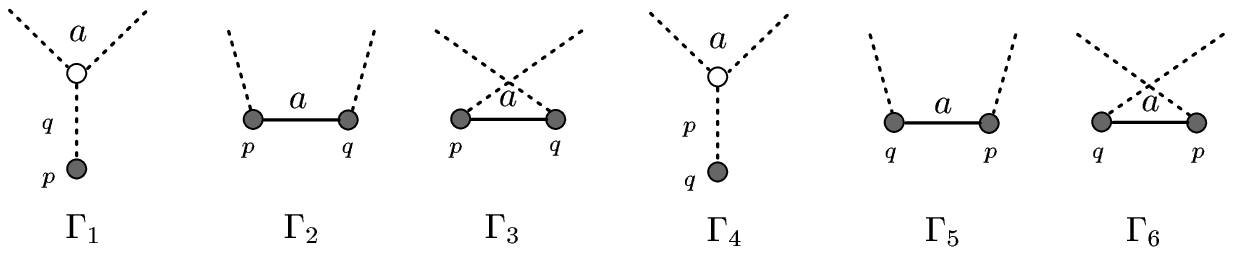}
\caption{The possible labelled graphs which give the same graph $\Gamma'$ after contractions of the middle edges.}\label{fig:stu_labeled}
\end{figure}
The six graphs are all possible ones which yield the same labelled graph $\Gamma'$ when the middle edges are contracted. The principal face contribution for $\Gamma_1$ with the middle $\theta$-edge, say $e$, collapsed is given by $\pm \int_{S^{n-1}}\omega_e\wedge I(\Gamma')=\pm I(\Gamma')$ while the contribution for $\Gamma_2$, $\Gamma_3$ with the middle $\eta$-edge, say $e'$, collapsed is given by $\pm I(\Gamma_2/e')$, $\pm I(\Gamma_3/e')=\pm I(\Gamma')$. The cases of $\Gamma_4$, $\Gamma_5$, $\Gamma_6$ are similar. The orientation of $\Sigma_A\cong S^{n-1}\times C_{\Gamma'}$ induced from $\mathrm{or}'(\Gamma_1)=\rho_{1}(a)\wedge\rho_{0}(p)\wedge\rho_0(q)\wedge O'$ (re-arranged in this form) is given by
\[  vol_{S^{n-1}}\wedge 
		i\Bigl(\frac{\partial}{\partial \rho_0(q)}\Bigr)
		i\Bigl(\frac{\partial}{\partial \rho_0(p)}\Bigr)
		i\Bigl(\frac{\partial}{\partial \rho_{1}(a)}\Bigr)\,\mathrm{or}'(\Gamma_1)
	=vol_{S^{n-1}}\wedge O'. \]
For other graphs $\Gamma_i$, we get the same $\mathrm{or}'(\Gamma_i)=\rho_{1}(a)\wedge\rho_{0}(p)\wedge\rho_0(q)\wedge O'$ and the induced orientation on $C_{\Gamma'}$ is again given by $O'$. Therefore we see that the terms 
$\sum_{i=1}^6[\Gamma_i](\pi_{\Gamma_i})^\partial_* \omega_{\Gamma_i}$ 
in the sum in (\ref{eq:dz}) restricted to the corresponding (principal) face of $C_{s+t}(\R^n)$ is of the form
\[ \bigl(\sum_{i=1}^6[\Gamma_i]\bigr)I(\Gamma')
	=2\bigl([\Gamma_1]+[\Gamma_2]+[\Gamma_3]\bigr)I(\Gamma'), \]
which vanishes by the STU relation $[\Gamma_1]+[\Gamma_2]+[\Gamma_3]=0$. 

The vanishing on other faces are shown in the rest of this section. Here we only give a guide to the rest of this section. The vanishing of the contributions of (3), the infinite faces, are shown by dimensional arguments (this has been shown in \cite[\S{5.8}]{S}). The vanishing of the contributions of (2), the hidden faces and when $n-j$ even the contribution of (4), anomalous faces, are discussed from the next subsection. In particular, through Lemmas~\ref{lem:hidden-even}, \ref{lem:anomaly-even}, \ref{lem:quasianomalous}. This will be the most complicated part in the proof. Finally when both $n$ and $j$ are even, we can not prove the vanishing on the anomalous faces (4). Fortunately, we can find the correction term as in the statement of Theorem~\ref{thm:closed} that kills the anomalous face contribution. It will be discussed in \S\ref{ss:anomaly_correction}.
\end{proof}

\subsection{Vanishing on hidden/anomalous faces, even codimension case}

 When the codimension is even and $\geq 2$, the following lemma immediately follows from lemmas given in \cite{S}, which is based on the codimension two case of \cite{R} (see also \cite{Wa}).
\begin{Lem}\label{lem:hidden-even}Suppose that the codimension is even and $\geq 2$.
Then the fiber integrals $(\pi_{\Gamma}^{\partial_A})_*\omega_\Gamma$, $A\subsetneq V(\Gamma)$, vanish.
\end{Lem}
Thus in the even codimension case the only contribution of ${\pi_\Gamma^\partial}_*\omega_\Gamma$ over non-principal faces is the contribution of the anomalous face. If moreover both $n$ and $j$ are odd, then the following lemma holds (see \cite[Proposition~5.17]{S}, \cite[Proposition~A.13]{Wa}).
\begin{Lem}\label{lem:anomaly-even}
If $n$ and $j$ with $n-j\geq 2$ are both odd, then the anomalous faces contribution vanishes, i.e., $dz_k=0$. Hence we have a well-defined cohomology class $[z_k]\in H^*(\emb{n}{j};\calA_k)$.
\end{Lem}
This shows Theorem~\ref{thm:closed} for $n, j$ odd case. For the case that both $n$ and $j$ are even see \S\ref{ss:anomaly_correction}.

\subsection{Vanishing on most of hidden/anomalous faces, odd codimension case}
Let $j,n$ be a pair of positive integers with codimension odd $\geq 3$. In this case almost all hidden faces contributions vanish (\cite[\S{5.7}]{S}), but we still need to prove the vanishings of contributions of other kinds of faces than those which do not contribute in the even codimension case, which correspond to the collapses of {\em admissible subgraphs}, to get a closed form on $\emb{n}{j}$. We say that a subgraph $\Gamma_A$ of an admissible graph $\Gamma$ is admissible if $\Gamma_A$ itself is admissible in the sense of Definition \ref{def:admissible} and if $|A|\geq 3$.

We will prove the following lemma in the rest of this subsection and the next subsection.

\begin{Lem}\label{lem:quasianomalous}
Suppose one of the following conditions holds:
\begin{itemize}
\item $n$ is odd and $j$ is even.
\item $n$ and $j$ satisfies the condition (1)-(b) or (1)-(c).
\end{itemize}
\noindent Then the fiber integrals $(\pi_\Gamma)^\partial_*\omega_\Gamma$ restricted to faces of $\partial C_{\Gamma}$ corresponding to the collapses of admissible subgraphs cancel each other in the sum $z_k$.
\end{Lem}

In the proof of Lemma~\ref{lem:quasianomalous} we will need the following lemma.
\begin{Lem}\label{lem:split}
For a subset $A\subset V(\Gamma)$, suppose that $\Gamma_A$ has an $\eta$-edge $e$ such that $\Gamma_A\setminus e$ is a disjoint union of two subgraphs $\Gamma_{A,1}$ and $\Gamma_{A,2}$ one of which has vertices at least two. Then $I(\Gamma)$ restricted to $\Sigma_A$ vanishes.
\end{Lem}
\begin{proof}
Let us consider the action of $\R_{>0}$ on $\Sigma_A$ given by dilations of points corresponding to vertices of $\Gamma_{A,2}$ around the intersection (point) of $\Gamma_{A,2}$ and $e$. The action of $\R_{>0}$ is free because $|A|\geq 3$. So we can consider the quotient $q: \Sigma_A\to \Sigma_A/\R_{>0}$ and it is easy to check that $\omega_{\Gamma_A}$ is basic with respect to $q$. The dimension of the fiber $\Sigma_A/\R_{>0}$ is strictly less than that of $\Sigma_A$. So the fiber integral vanishes by a dimensional reason.
\end{proof}

\begin{proof}[Proof of Lemma~\ref{lem:quasianomalous} (partial)]
Suppose that $|A|\geq 3$ and that the subgraph $\Gamma_A\subsetneq \Gamma$ is admissible. 


Let us first suppose that $\Gamma_A$ is a tree.
If moreover $\Gamma_A$ has an $\eta$-edge, then the vanishing follows from Lemma~\ref{lem:split} above. 

If $\Gamma_A$ is a $Y$-shaped admissible graph with only $\theta$-edges, then the vanishing of the integral is implied by the Y relation. In this case, six labelled graphs cancel each other. 

If $\Gamma_A$ is a tree with only $\theta$-edges and with at least two e-vertices, then $\Gamma_A$ has a subgraph $\Gamma_I$ as depicted in Figure~\ref{fig:GI} (all the i-vertices in the figure are univalent in $\Gamma_A$).
\begin{figure}
\includegraphics{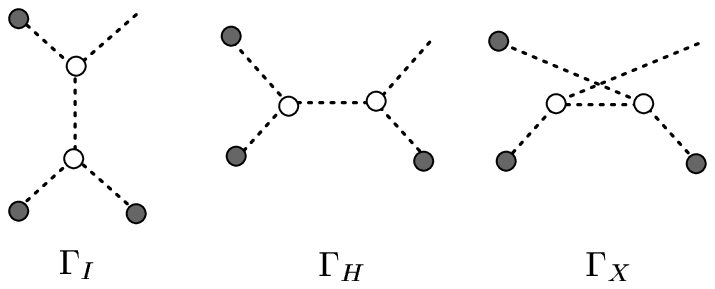}
\caption{The three subgraphs which cancel each other.}\label{fig:GI}
\end{figure}
There are other possibilities for $\Gamma$'s which agree with $\Gamma$ except for the subgraph $\Gamma_I$ replaced by $\Gamma_H$ or $\Gamma_X$ as depicted in Figure~\ref{fig:GI} with labels as given in the relation in Figure~\ref{fig:relations_odd}. Let us denote these graphs by $\Gamma'$, $\Gamma''$. It is easy to check that the integrals of $\Gamma$, $\Gamma'$, $\Gamma''$ coincide on the face $\Sigma_{A}$. Hence in the labelled graph expression of $z_k$ we see that 
\[ [\Gamma](\pi_\Gamma^{\partial_A})_*\omega_{\Gamma_A}
	+[\Gamma'](\pi_{\Gamma'}^{\partial_A})_*\omega_{\Gamma_A'}
	+[\Gamma''](\pi_{\Gamma''}^{\partial_A})_*\omega_{\Gamma_A''}
	=([\Gamma]+[\Gamma']+[\Gamma''])(\pi_\Gamma^{\partial_A})_*\omega_{\Gamma_A}=0 \]
by the IHX relation.


Next we suppose that $\Gamma_A$ is not a tree. In this case either 
\begin{itemize}
\item $\Gamma_A=\Gamma_{A,1}\cup e\cup \Gamma_{A,2}$ where $e$ is an $\eta$-edge, $\Gamma_{A,1}$ is a tree, $\Gamma_{A,2}$ has a loop, and $\Gamma_{A,1}\cap\Gamma_{A,2}=\emptyset$, or
\item $\Gamma_A$ has a part as in Figure~\ref{eq:strict_admissible}. 

\begin{figure}
\begin{equation*}
\includegraphics{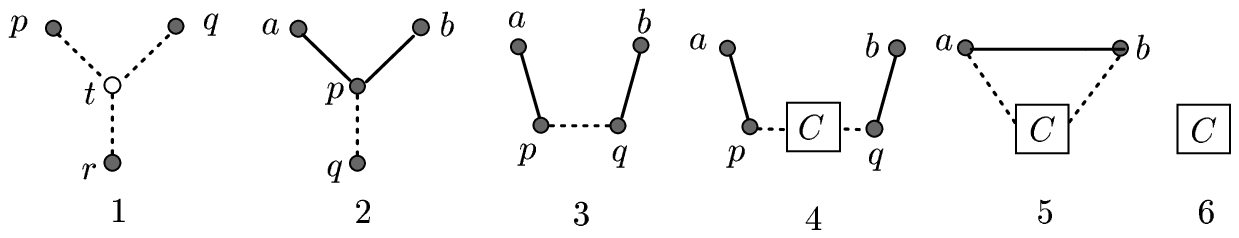}
\end{equation*}
\caption{$C$ is a component with only $\theta$-edges.  In each graph there are no other edges incident to $p$ and $q$ than those shown there, and $r$, $a$, $b$ are not univalent.}\label{eq:strict_admissible}
\end{figure}
\end{itemize}

Now we show the vanishing for each of these cases.

\begin{enumerate}
\item If $\Gamma_A=\Gamma_{A,1}\cup e\cup \Gamma_{A,2}$ as in the first case, then the vanishing of the integral follows again from Lemma~\ref{lem:split} above.

\item If $\Gamma_A$ has a subgraph of type 1 in Figure~\ref{eq:strict_admissible}, then it must be that one or two $\eta$-edges share the vertex $r$. If it is just one, then the vanishing follows from Lemma~\ref{lem:split} above. If it is just two, then let $(r,a)$ and $(r,b)$ be the two $\eta$-edges. Consider the automorphism $g:\hat{B}_A\to \hat{B}_A$ given as follows:
\[\begin{split}
	 g:(f;x_a,x_b,x_p,x_q,x_r,x_t,\ldots)
	&\mapsto (f;x_a,x_b,x_p+(x_a+x_b-2x_r),x_q+(x_a+x_b-2x_r),\\
	&\hspace{11mm}x_a+x_b-x_r,x_t+f(x_a+x_b-2x_r),\ldots).
	\end{split} \]
This can be realized by a central symmetry of $x_r$ around the center of $x_ax_b$ ($x_r\mapsto x_a+x_b-x_r$) followed by translations of $x_p,x_q,x_t$ by the difference $(x_a+x_b-x_r)-x_r$. If $n$ even $j$ odd, $g$ reverses the orientation of the fiber and preserves the sign of $\hat\omega_{\Gamma_A}$, i.e., $g^*\hat\omega_{\Gamma_A}=\hat\omega_{\Gamma_A}$. If $n$ odd $j$ even, then $g$ preserves the orientation of the fiber and reverses the sign of $\hat\omega_{\Gamma_A}$. Hence the integral vanishes.

\item If $\Gamma_A$ has a subgraph of type 2 or 3 in Figure~\ref{eq:strict_admissible}, consider the automorphism $g:\hat{B}_A\to \hat{B}_A$ given by
\[ g: (f;x_a,x_b,x_p,x_q,\ldots)\mapsto (f;x_a,x_b,x_a+x_b-x_q,x_a+x_b-x_p,\ldots). \]
(This symmetry has been considered in \cite[Lemma~6.5.5]{R}.) When $n$ odd $j$ even, $g$ preserves the orientation of the fiber and reverses the sign of the integrand form. When $n$ even $j$ odd, $g$ reverses the orientation of the fiber and preserves the sign of the integrand form. Hence in any case the integral vanishes.

\item If $\Gamma_A$ has a subgraph of type 4 in Figure~\ref{eq:strict_admissible}, consider the symmetry of $\hat{B}_A$ given by the composition of the following symmetries:
\begin{enumerate}
\item Central symmetry of the subgraph between $a$ and $b$ around the point $\frac{x_a+x_b}{2}$. Write $p'$ and $q'$ the images of $p$ and $q$ respectively.
\item Central symmetry of the inverted subgraph between $p'$ and $q'$ around the point $\frac{x_{p'}+x_{q'}}{2}$.
\end{enumerate}
One can check the vanishing of the integral as in the type 3 case.

\item The case when $\Gamma_A$ has a subgraph of type 5 in Figure~\ref{eq:strict_admissible} or of type 6 will be separately discussed in the next subsection.\qedhere
\end{enumerate}
\end{proof}

\subsection{Vanishing for type 5 or 6 subgraphs, odd codimension case}

We continue to study the odd codimension case. Now we consider in particular the case where an admissible subgraph $\Gamma_A$ does not have an $\eta$-edge (type 6), or has just one $\eta$-edge (type 5, see Figure~\ref{eq:strict_admissible}). We will call such a $\Gamma_A$ an {\it special subgraph}. We show that a sum of special graphs contributions cancel each other in some sense generalizing the cancelling argument of the principal faces contributions, given in \S\ref{s:outline-proof}.

\subsubsection{Local description of $z_k$}
If $\Gamma_A$ is special, then we may assume that it consists of a type (a) path (see Figure~\ref{fig:paths}) with some hairs replaced by $Y$-shaped graphs (as the graphs in Example 2 below) and at most one $\eta$-edge. This is because special graphs with more complicated trees consisting only of $\theta$-edges cancel each other as shown in Figure~\ref{fig:GI}. In the following we assume that $\Gamma_A$ is special of order $\ell$. 

We have seen that the configuration space integral $(\pi_\Gamma)_*^\partial\omega_\Gamma$ restricted to the face $\Sigma_A$ is expressed as
\begin{equation}\label{eq:integral_sigma}
 \int_{C_{\Gamma/\Gamma_A}}
		D_A^*\rho_{A*}\hat{\omega}_{\Gamma_A}\wedge \omega_{\Gamma/\Gamma_A} 
\end{equation}
(See (\ref{eq:Sigma_A})). We would like to show that a linear combination of the integrals of this form vanishes. We claim that a cancel occurs among the terms (\ref{eq:integral_sigma}) for pairs $(\Gamma', \Gamma_B')$ such that $\Gamma'\in\tilde\calG_k$, $\Gamma'_B$ admissible subgraph of $\Gamma'$ and $\Gamma'/\Gamma'_B=\Gamma/\Gamma_A$ for a fixed pair $(\Gamma,\Gamma_A)$. 

To see this we fix the data $Q=(\Gamma^Q, v, \ell)$ where
\begin{enumerate}
\item $\Gamma^Q:=\Gamma/\Gamma_A$ for some admissible pair $\Gamma_A\subset \Gamma$, $\Gamma\in\tilde\calG_k$, equipped with a suitable label and with one vertex $v\in V(\Gamma^Q)$ distinguished as the point where $\Gamma_A$ is collapsed,
\item $\ell=\mathrm{ord}(\Gamma_A)=|A|/2$.
\end{enumerate}
Note that there may be several possibilities for $\Gamma$ of order $k$ and its admissible subgraph $\Gamma_A$ of order $\ell$ that yield the same triple as $Q$. We consider all such order $\ell$ admissible subgraphs of graphs in $\tilde\calG_k$ that yield the same triple as $Q$. We denote by $\tilde\calG_\ell(Q)$ the set of all such admissible subgraphs and let $\calG_\ell(Q)=\R\tilde\calG_\ell(Q)/(\Gamma,-\mathrm{or})=-(\Gamma,\mathrm{or})$. Note that graphs in $\calG_\ell(Q)$ are subgraphs. So we forget external structure. Then consider the following $\calG_\ell(Q)$-linear combination of the integrands $D_A^*\rho_{A*}\hat{\omega}_{\Gamma_A}$ for such graphs: 
\[
 z_\ell'(Q):=\sum_{{\Gamma_A}\atop{\mathrm{labelled}}}
		\Gamma_A\otimes D_A^*\rho_{A*}\hat{\omega}_{\Gamma_A}
		\in \calG_\ell(Q)\otimes\Omega_{DR}^*(C_{\Gamma^Q})
\]
where the sum is taken over admissible subgraphs in $\tilde\calG_\ell(Q)$.

\begin{figure}
\fig{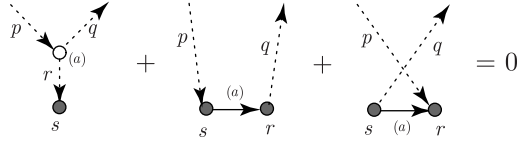}
\caption{STU' relation}\label{fig:stu2}
\end{figure}
Let $\calA_\ell(Q)$ be the space of $\Gamma_A$'s in $\tilde{\calG}_\ell(Q)$ labelled oriented, quotiented by the ``labelled versions" of the IHX, ST2, STU, Y, L and the STU' relation (Figure~\ref{fig:stu2}, the ST relation and the label change relation are excluded). Namely, the 2- or 3-term relations given in Figure~\ref{fig:relations_odd} are the ones obtained from the 4- or 6-term relations by modding out the label changes. The labelled relations we consider here is the 4- or 6-term relations. Now we define the following maps:
\begin{enumerate}
\item The map $i_Q:\calG_\ell(Q)\to \calG_k$ is defined for $\Gamma_A\in\tilde\calG_\ell(Q)$ by the sum of all possible admissible replacements of the vertex $v$ of $\Gamma^Q$ with $\Gamma_A$.
\item The map $m_2:\calG_k\to \calG_k$ is defined for $\Gamma\in\tilde\calG_k$ by $m_2(\Gamma)=2^p\Gamma$ where $p$ is the number of univalent vertices of $\Gamma$. This will be necessary in order that STU' relations are mapped to ST relations.
\end{enumerate}
Then  by comparing the defining relations for $\calA_\ell(Q)$ and $\calA_k$ we have the following Lemma.
\begin{Lem}\label{lem:mi}
The map $m_2\circ i_Q:\calG_\ell(Q)\to \calG_k$ descends to a well-defined map $\overline{i}_Q:\calA_\ell(Q)\to \calA_k$.
\end{Lem}
Lemma~\ref{lem:mi} shows that if we define
\begin{equation}\label{eq:z(A)}
 z_\ell(Q):=([\cdot]\otimes 1)(z_\ell'(Q)) \in \calA_\ell(Q)\otimes \Omega_{DR}^*(C_{\Gamma^Q}) 
\end{equation}
then $\int_{C_{\Gamma^Q}}(\overline{i}_{Q}\otimes 1)(z_\ell(Q))\wedge \omega_{\Gamma^Q}$ is a constant multiple of a partial sum in the formula (\ref{eq:dz}) of $dz_k$ restricted to $\Sigma_A$'s and $dz_k$ restricted to $\Sigma_A$ is a sum of such terms. So it is enough for our purpose to show that $z_\ell(Q)=0$ for any $Q$. Note that from the discussion above, we see that only the special graph terms survive in $z_\ell(Q)$.

\subsubsection{Decomposition to units}

To study $z_\ell(Q)$, we decompose the set of special graphs into small pieces. It is observed that if a special subgraph $\Gamma_A$ of $\Gamma$
\begin{enumerate}
\item does not have an $\eta$-edge, then by the IHX relation it is expanded in a sum of $\ell$-wheels in $\calA_\ell(Q)$ where an {\it $\ell$-wheel} is a labelled graph whose underlying graph is shown in Figure~\ref{fig:std_lab_wheel} (with possibly different labels from that of the figure).
\item has an $\eta$-edge, then by the ST2/STU relation there is another labelled special (sub)graph $\Gamma_A'$ (of $\Gamma'$), which differs from $\Gamma_A$ only by an orientation preserving label change, so that $\Gamma_A+\Gamma_A'$ is equivalent in $\calA_\ell(Q)$ to a sum of graphs without $\eta$-edges. Then $\Gamma_A+\Gamma_A'$ is expanded in $\calA_\ell(Q)$ in a sum of $\ell$-wheels.
\end{enumerate}
This observation suggests a decomposition of the set $\tilde\calG_\ell(Q)$ of special graphs into pieces, which we will call {\it units}. Namely by a {\it unit} we mean a single graph $\Gamma_A$ in the case (1) above, or a pair of graphs $(\Gamma_A,\Gamma_A')$ as above in the case (2). Then by definition a sum of terms in a single unit is equivalent in $\calA_\ell(Q)$ to a sum of $\ell$-wheels.

Since a special subgraph has at most one $\eta$-edge, no two different units overlaps. Hence the set $\tilde\calG_\ell(Q)$ is decomposed into disjoint units. Below we shall prove the cancelling between one or two units, which will conclude $z_\ell(Q)=0$. 

\subsubsection{Cyclic permutation of a label on $\Gamma_A$}
Now let us assume that $n$ is odd and $j$ is even and that $\Gamma_A$ is special. The case where $n$ is even and $j$ is odd will be discussed later in page \pageref{pg:even-odd}. We can first see that the hidden face contribution of $\Sigma_A$ with $\Gamma_A$ being odd order vanishes. This is because the central symmetry in $\R^n$ of the local configuration space with respect to one of points lying on the $j$-dimensional plane $f(\R^j)$ (as in the proof of \cite[Proposition~A.13]{Wa}) reverses the orientation of the fiber and preserves the sign of the integrand form.

The same argument does not work when the special subgraph $\Gamma_A$ is of even order. Instead we prove the vanishing for terms of even order subgraphs by considering a cyclic permutation symmetry acting simultaneously on all graphs in a unit. A `cyclic permutation' of a label on $\Gamma_A$ is defined as follows. As in Definition~\ref{def:labelled_graph} one can also define $S(\Gamma_A)$ and $T(\Gamma_A)$ for $\Gamma_A$, namely, $S(\Gamma_A)=V_e(\Gamma_A)\sqcup E_\eta(\Gamma_A)$, $T(\Gamma_A)=V_i(\Gamma_A)\sqcup E_\theta(\Gamma_A)$. Recall that $S$-labelled (resp. $T$-labelled) objects are of odd degree (resp. even degree). We consider that a label on $\Gamma_A$ is given by numberings on the sets $S(\Gamma_A)$ and $T(\Gamma_A)$. As for graphs in $\tilde{\calG}_k$, a label on $\Gamma_A$ together with a choice of an orientation of each $\theta$-edge determines an orientation of $\Gamma_A$. 

\begin{figure}
\includegraphics{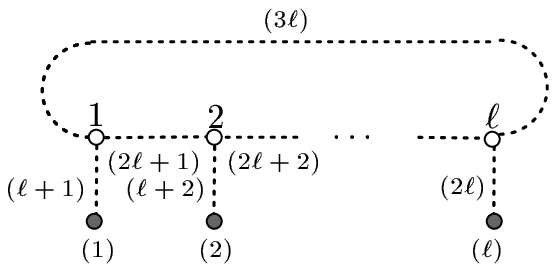}
\caption{Standard labelling on a $\ell$-wheel.}\label{fig:std_lab_wheel}
\end{figure}

\begin{figure}
\fig{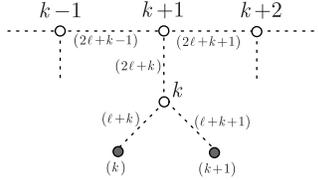}
\caption{Standard labelling on a non-wheel special graph without $\eta$-edges}\label{fig:std_order}
\end{figure}

There is a natural choice of a cyclic ordering on the set $S(\Gamma_A)$ given as follows. If $\Gamma_A$ is a labelled $\ell$-wheel, then $S(\Gamma_A)=V_e(\Gamma_A)$ and the natural cyclic ordering is defined by the standard labelling given in Figure~\ref{fig:std_lab_wheel}. For non-wheel special subgraphs without $\eta$-edges, the standard labelling is given as in Figure~\ref{fig:std_order}. For non-wheel special subgraphs with an $\eta$-edge, namely for type 5 graphs of Figure~\ref{eq:strict_admissible}, natural cyclic orderings are canonically induced from those of an $\ell$-wheel: in the STU relation, for example, if one of the three terms in the relation is given a $S$-label then the $S$-labels of the others are canonically determined so that these are compatible with the graph orientations that are consistent with the STU relation. See Figure~\ref{fig:stu_labeled}.

The natural cyclic ordering defines a set automorphism
\[ 
	\sigma:S(\Gamma_A)\to S(\Gamma_A)
\]
given by taking the next element with respect to the (increasing) order. This turns $\Gamma_A$ into another labelled graph by changing an $S$-label $P$ into $\sigma^{-1}(P)$. 
If we change the label, the automorphism $\sigma$ changes the label of $\Gamma_A$ and so may change the sign of the integral $D_A^*\rho_{A*}\hat{\omega}_{\Gamma_A}$ (with respect to the corresponding automorphism of the configuration space). More precisely, according to the definition of the integral in \S\ref{subsec:integral}, a cyclic permutation of the $S$-label induced by $\sigma$ acts on the fiber integral as $-1$ because the sign of an even cyclic permutation (of odd elements) is $-1$.

\begin{proof}[Proof of Lemma~\ref{lem:quasianomalous} (continued), $n$ odd, $j$ even, $\ell$ even case]

As we have observed, we need only to prove the cancelling of the integrals restricted to the faces corresponding to collapses of special subgraphs. Suppose, for simplicity, that the set $S(\Gamma_A)$ is labelled by $\{1,2,\ldots,\ell\}$ so that $1<2<\cdots<\ell<1$ in the natural cyclic ordering given above. The other cases can be treated separately and analogously. Let $\tilde\calG_\ell^\mathrm{std}(\Gamma_A)$ be the set of labelled special subgraphs in $\tilde\calG_\ell(Q)$ with isomorphic underlying edge-oriented unlabelled graph as $\Gamma_A$, and with the labelling on $S(\Gamma_A)$ satisfying the simplicity assumption above. 

Now take a unit $u(\Gamma_A)$ and write as $u(\Gamma_A)=\Gamma_A^*\in \tilde{\calG}_\ell^\mathrm{std}(\Gamma_A)$ if $|u(\Gamma_A)|=1$, or as $u(\Gamma_A)=(\Gamma_A^*,\Gamma_A^{**})\in \tilde{\calG}_\ell^\mathrm{std}(\Gamma_A)^{\times 2}$ if $|u(\Gamma_A)|=2$, and expand $\Gamma_A^*$ or $\Gamma_A^*+\Gamma_A^{**}$ in a sum of $\ell$-wheels in $\calA_\ell(Q)$: $\Gamma^*_{u,1}+\Gamma^*_{u,2}+\cdots+\Gamma^*_{u, N}$ ($\Gamma_{u,i}^*$: $\ell$-wheel). This expansion is unique up to permutations of suffixes $i=1,\ldots,N$, and the correspondence 
\begin{equation}\label{eq:expand_u1}
 \mbox{(a labelling $\rho$ on $u(\Gamma_A)$)}\mapsto (\Gamma_{u,1}^*(\rho),\Gamma_{u,2}^*(\rho),\ldots,\Gamma_{u,N}^*(\rho)) 
\end{equation}
determines (non-uniquely) a matrix $M$ (each labelling corresponds to a row of $M$) where $\Gamma_{u,i}^*(\rho)$ is $\Gamma_{u,i}^*$ with the induced labelling. We view $M$ as a multiset consisting of labelled oriented wheels.

For each fixed $\Gamma^*_{u,i}(\rho)$ in (\ref{eq:expand_u1}), there is a non-identity permutation 
\[ \tau:T(\Gamma_A^*)\to T(\Gamma_A^*) \]
acting on the $T$-label(s) of graph(s) of $u(\Gamma_A)$ defined so that the $\ell$-wheel expansion of $\tau\sigma u(\Gamma_A)$ in the labelling $\rho$: $\tau\sigma \Gamma_{u,1}^{*}(\rho)+\tau\sigma \Gamma_{u,2}^{*}(\rho)+\cdots+\tau\sigma \Gamma_{u,N}^{*}(\rho)$ has a term $\tau\sigma\Gamma_{u,j}^*(\rho)$ with
\begin{equation}\label{eq:G=G}
 [\Gamma_{u,i}^*(\rho)]=[\tau\sigma\Gamma_{u,j}^{*}(\rho)]=[\Gamma_{u,j}^{*}(\tau\sigma\rho)]. 
\end{equation}
Note that $\tau$ is uniquely determined by $\Gamma_{u,i}^*(\rho)$: the labelled graph $\sigma\Gamma_{u,i}^*(\rho)$ is isomorphic to the labelled graph obtained from $\Gamma_{u,i}^*(\rho)$ by a permutation $\varphi$ on $T(\Gamma_{u,i}^*(\rho))$ (keeping $S$-labels fixed). Then $\tau:T(\Gamma_A^*)\to T(\Gamma_A^*)$ is given by $\tau(x)=\varphi^{-1}(x)$ where $T(\Gamma_A^*)$ is naturally identified with $T(\Gamma_{u,i}^*(\rho))$ by the labels.

Now in the $\ell$-wheel expansion of the sum $z_\ell(Q)$ we see that the terms for $\Gamma_{u,i}^{*}(\rho)$ and $\tau\sigma\Gamma_{u,j}^*(\rho)$ cancel each other, i.e.,
\[ 
	[\Gamma_{u,i}^*(\rho)]\otimes D_A^*\rho_{A*}\hat{\omega}_{\Gamma_{A}^*}
	+ [\tau\sigma\Gamma_{u,j}^*(\rho)]\otimes D_A^*\rho_{A*}\hat{\omega}_{\tau\sigma \Gamma_{A}^*}
	=([\Gamma_{u,i}^*(\rho)]-[\tau\sigma\Gamma_{u,j}^*(\rho)])\otimes D_A^*\rho_{A*}\hat{\omega}_{\Gamma_{A}^*}=0
	 \]
by (\ref{eq:G=G}) and by the fact that $\sigma$ only changes the sign of the integral and that $\tau$ does not change the integral (though they may change the coefficient graph). More generally, the mapping $\Gamma_{u,j}^*(\rho)\mapsto \tau\sigma\Gamma_{u,j}^*(\rho)$ ($\tau$ depends on $\Gamma_{u,j}^*(\rho)$) induces an automorphism on the multiset $M$ without fixed point. Hence the cancelling pairs are mutually disjoint and all terms in $M$ cancel with each other. Note that the sum $z_\ell(Q)$ is over the rows of $M$ (one row for one term) for each unit $u(\Gamma_A)$.
\end{proof}

\begin{proof}[Example~1]
Let us see some typical examples for the cancellation. We assume that $n$ odd, $j$ even. First by the STU/ST2 relation, we have the following identities
\begin{equation}\label{eq:ex_cancel}
\raisebox{-0.45\height}{\includegraphics{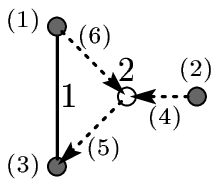}}
 +\raisebox{-0.45\height}{\includegraphics{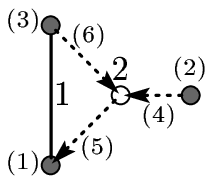}}
=\raisebox{-0.35\height}{\includegraphics{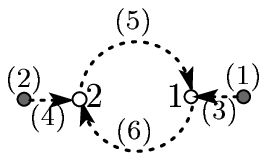}}
=\raisebox{-0.45\height}{\includegraphics{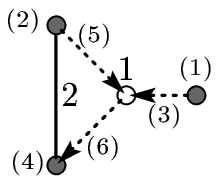}}
 +\raisebox{-0.45\height}{\includegraphics{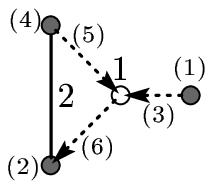}}
\end{equation}
in $\calA_2(Q)$. Let $u_1=(\Gamma_1,\Gamma_1')$ be the unit consisting of the first two graphs of (\ref{eq:ex_cancel}) and let $u_2=(\Gamma_2,\Gamma_2')$ be that of the last two graphs. Then it holds that 
$ u_2=\tau\sigma u_1 $
where $\sigma$ is the cyclic permutation acting on the set $S=\{1,2\}$, $ \tau=(1\ 2)(3\ 4)(5\ 6)$, 
$T=\{1,2,3,4,5,6\}$, and that $(\tau\sigma)^2=\mathrm{id}$. Then we see that
\[ \begin{split}
	&[\Gamma_1]\otimes D_A^*\rho_{A*}\hat\omega_{\Gamma_1}
	+[\Gamma_1']\otimes D_A^*\rho_{A*}\hat\omega_{\Gamma_1'}
	+[\Gamma_2]\otimes D_A^*\rho_{A*}\hat\omega_{\Gamma_2}
	+[\Gamma_2']\otimes D_A^*\rho_{A*}\hat\omega_{\Gamma_2'}\\
	&=([\Gamma_1]+[\Gamma_1']-[\Gamma_2]-[\Gamma_2'])\otimes  D_A^*\rho_{A*}\hat\omega_{\Gamma_1}=0
	\end{split} \]
by the relation (\ref{eq:ex_cancel}). The contribution of any other special graph of order 2 with one $\eta$-edge is cancelled by the same argument.
\end{proof}

\begin{proof}[Example~2]
Assume that $n$ odd, $j$ even again. Consider the special graphs (units)
\[
 X(p,q;r,s):=\raisebox{-0.45\height}{\includegraphics{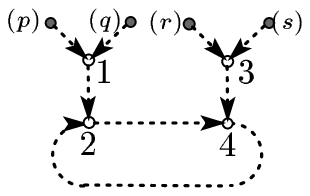}},\quad
 Y(p,q;r,s):=\raisebox{-0.45\height}{\includegraphics{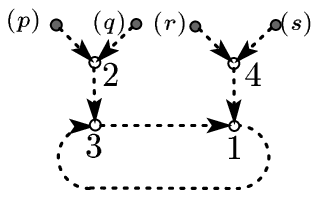}}
\]
where $\{p,q,r,s\}$ is a permutation of $\{1,2,3,4\}$. $X(p,q;r,s)$ and $Y(p,q;r,s)$ are related to each other by $\sigma$. One may fix a standard way of labelling on edges of $X$'s and $Y$'s from $p,q,r,s$. So we fix one such. The cases of other choices can be discussed similarly. Let
\[ W(p,q,r,s):=\raisebox{-0.45\height}{\includegraphics{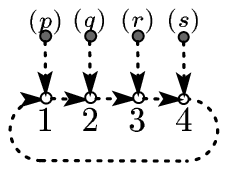}}. \]
Then by the IHX relation we have
\begin{equation}\label{eq:list_wheels}
 \begin{split}
	X(1,2;3,4)&=W(1,2,3,4)+W(2,1,3,4)+W(1,2,4,3)+W(2,1,4,3)\\
	X(3,4;1,2)&=W(3,4,1,2)+W(3,4,2,1)+W(4,3,1,2)+W(4,3,2,1)\\
	X(1,3;2,4)&=W(1,3,2,4)+W(3,1,2,4)+W(1,3,4,2)+W(3,1,4,2)\\
	X(2,4;1,3)&=W(2,4,1,3)+W(2,4,3,1)+W(4,2,1,3)+W(4,2,3,1)\\
	X(1,4;2,3)&=W(1,4,2,3)+W(4,1,2,3)+W(1,4,3,2)+W(4,1,3,2)\\
	X(2,3;1,4)&=W(2,3,1,4)+W(2,3,4,1)+W(3,2,1,4)+W(3,2,4,1)\\
	Y(1,2;3,4)&=W(4,1,2,3)+W(4,2,1,3)+W(3,1,2,4)+W(3,2,1,4)\raisebox{.5\height}{$\phantom{\int}$}\\
	Y(3,4;1,2)&=W(2,3,4,1)+W(1,3,4,2)+W(2,4,3,1)+W(1,4,3,2)\\
	Y(1,3;2,4)&=W(4,1,3,2)+W(4,3,1,2)+W(2,1,3,4)+W(2,3,1,4)\\
	Y(2,4;1,3)&=W(3,2,4,1)+W(1,2,4,3)+W(3,4,2,1)+W(1,4,2,3)\\
	Y(1,4;2,3)&=W(3,1,4,2)+W(3,4,1,2)+W(2,1,4,3)+W(2,4,1,3)\\
	Y(2,3;1,4)&=W(4,2,3,1)+W(1,2,3,4)+W(4,3,2,1)+W(1,3,2,4)\\
	\end{split}
 \end{equation}
in $\calA_4(Q)$. Here $\sigma$ maps $X(p,q;r,s)$ to $Y(p,q;r,s)$ and $\sigma$ acts on wheels. For example, $\sigma$ maps $W(1,2,3,4)$ to $W(4,1,2,3)$ and for this term $\tau=(1\ 2\ 3\ 4)$. In this case $\tau\sigma X(1,2;3,4)=Y(2,3;4,1)=Y(2,3;1,4)$. Indeed the expansion of $Y(2,3;1,4)$ includes $W(1,2,3,4)$ too. Noting that the integrals for $X(p,q;r,s)$ are all equal, say to $\alpha$, and that the integrals for $Y(p,q;r,s)$ are all equal to $-\alpha$ by definition of integral in \S\ref{subsec:integral}, it follows easily by using (\ref{eq:list_wheels}) that
\[ \sum_{(p,q;r,s)}\Bigl(
	[X(p,q;r,s)]\otimes D_A^*\rho_{A*}\hat{\omega}_{X(p,q;r,s)}
	+[Y(p,q;r,s)]\otimes D_A^*\rho_{A*}\hat{\omega}_{Y(p,q;r,s)}\Bigr)=0.\qedhere \]
\end{proof}

\begin{proof}[Proof of Lemma~\ref{lem:quasianomalous} (continued), the case $n$ even, $j$ odd]\label{pg:even-odd}
We consider the following cases as given in the statement of Theorem~\ref{thm:closed}: (1)-(b) $k\leq 4$,
(1)-(c) $j=3$ and $n\geq 12$.

In the case (1)-(b), the vanishing of the contributions of $\Gamma_A$'s of even order can be shown similarly as in the case $n$ odd, $j$ even, $\ell$ odd by using the central symmetry around a univalent vertex. The vanishing of $\Gamma_A$'s of order 3 can be shown by replacing the cyclic permutation in the discussion above with the symmetry that reverses a 3-wheel around an axis. Note that the same argument does not work for $\ell\equiv 1\mbox{ mod }4$. So $(\ell\leq)k\leq 4$ is necessary.

However, in the special case as in (1)-(c), the vanishing can be proved for all $\ell$. The case $\ell=3$ has been done already. For $\Gamma_A$'s of order $\ell$ with $\ell\geq 5$, we have that $\deg \rho_{A*}\hat{\omega}_{\Gamma_A}=\ell(n-5)+4\geq 5n-21$. But when $n\geq 12$, we have that $5n-21 > \dim\calI_3(\R^n)=3n$. Therefore $\rho_{A*}\hat\omega_{\Gamma_A}=0$ by a dimensional reason.
\end{proof}

We have shown Lemma~\ref{lem:quasianomalous} so far and hence we have the following
\begin{Prop}\label{prop:dz=anomaly}
Suppose that $n, j, k$ satisfy one of the conditions in the statement of Theorem~\ref{thm:closed}. Then the exterior derivative of $z_k$ is rewritten as
\[ dz_k
	=\left\{
	\begin{array}{ll}
	\displaystyle\frac{1}{k_S!k_T!}\sum_{{\Gamma}\atop{\mathrm{labelled}}}[\Gamma]\otimes J\int_{C_1(\R^j)}D_{V(\Gamma)}^*\rho_{V(\Gamma)*}\hat{\omega}_\Gamma & \mbox{$n,j$: even}\\
	0 & \mbox{otherwise}
	\end{array}\right. \]
where $\int_{C_1 (\R^j)}$ denotes the integration along the fiber.\qed
\end{Prop}
This completes the proof of Theorem~\ref{thm:closed}(1).

\subsection{The anomalous face correction term}\label{ss:anomaly_correction}
In the rest of this section we let $A=V(\Gamma)$. As was observed in \S\ref{s:outline-proof} we know that the integral $I(\Gamma)$ restricted to the anomalous face $\Sigma_{A}$ can be written as the integral along $C_1(\R^j)$ of the differential form
\begin{equation}\label{eq:DArA}
 D_A^*{\rho_A}_*\hat\omega_{\Gamma} \in \Omega_\mathrm{DR}^{(n-j-2)k+j+1}(C_1(\R^j)\times \emb{n}{j}). 
\end{equation}

Now we would like to find an $(n-j-2)k+j$ form $\beta_{\Gamma}$ on $C_1(\R^j)\times \femb{n}{j}$ so that
\begin{equation}\label{eq:dbeta}
 \sum_{{\Gamma}\atop{\mathrm{labelled}}} [\Gamma ]\otimes d\int_{C_1(\R^j)}\beta_\Gamma
	=\sum_{{\Gamma}\atop{\mathrm{labelled}}} [\Gamma ]\otimes Jr^*\int_{C_1(\R^j)}D_A(\varphi)^* {\rho_A}_* \hat{\omega}_\Gamma.
\end{equation}
If such a $\beta_\Gamma$ is found, and if we set
\[ \Theta_k:=\frac{1}{k_S!k_T!}\sum_{{{\Gamma}\atop{\mathrm{labelled}}}}
	[\Gamma]\otimes\int_{C_1(\R^j)}\beta_\Gamma
	\in \calA_k\otimes\Omega_{DR}^{(n-j-2)k}(\femb{n}{j}), \]
then by Proposition \ref{prop:dz=anomaly}, the form $\hat{z}_k$ defined in (\ref{eq:hatz}) gives a closed $(n-j-2)k$-form on $\femb{n}{j}$, as desired in Theorem~\ref{thm:closed}(2) and completes the proof of Theorem~\ref{thm:closed}(2).


Recall that $\femb{n}{j}$ is the space of families $\widetilde{\varphi}=\{ \varphi_t \}$ of immersions
$\varphi_t:\R^j\to \R^n$, $t\in[0,1]$ such that $\varphi_0=\iota$ and $\varphi_1\in\emb{n}{j}$. We define a map
\[ \widetilde{D}_A:[0,1]\times C_1(\R^j)\times\femb{n}{j} \to \calI_j(\R^n) \]
by $\widetilde{D}_A(t,x,\widetilde{\varphi}=\{\varphi_t\})=D\varphi_t(x)$. Note that $D\varphi:T\R^j\to T\R^n$ is the differential of $\varphi$, which is linear injective when $\varphi$ is an immersion. $\widetilde{D}_A$ restricts on $\{0,1\}\times C_1(\R^j)\times\femb{n}{j}$ to $D_A(\iota)\circ (\id\times r)$ and $D_A(\varphi)\circ (\id\times r)$.

Then, put
\[ \beta_{\Gamma}:=-\mathrm{pr}_{23*}\widetilde{D}_{A}^* \rho_{A*}\hat\omega_{\Gamma}
	\in \Omega_\mathrm{DR}^{(n-j-2)k+j}(C_1(\R^j)\times\femb{n}{j}) \]
where $\mathrm{pr}_{23}:[0,1]\times C_1(\R^j)\times \femb{n}{j}\to C_1(\R^j)\times\femb{n}{j}$ is the projection.
\begin{Lem}
(\ref{eq:dbeta}) holds.
\end{Lem}
\begin{proof}
We use the generalized Stokes theorem (\ref{eq:stokes}); suppose $\deg \mathrm{pr}_{23*}\widetilde{D}_A^*\rho_{A*}\hat\omega_{\Gamma}=a$. Then we have
\begin{equation}\label{eq:dbeta_proof}
\begin{split}
		\sum_{\Gamma} [\Gamma ]\otimes
			d\beta_\Gamma
		&= - \sum_{\Gamma} [\Gamma ]\otimes
			d\,\mathrm{pr}_{23*}\widetilde{D}_A^*\rho_{A*}\hat\omega_{\Gamma}\\
		&= - \sum_{\Gamma} [\Gamma ]\otimes
			\left[ \mathrm{pr}_{23*}(d\,\widetilde{D}_A^*\rho_{A*}\hat\omega_{\Gamma})
			+ (-1)^{a+1} \mathrm{pr}_{23*}^\partial (\widetilde{D}_A^*\rho_{A*}\hat\omega_{\Gamma})
			\right]\\
		&= (-1)^a \sum_{\Gamma} [\Gamma ]\otimes
			\Bigl[(\id\times r)^*D_A(\iota)^*\rho_{A*}\hat\omega_{\Gamma}-(\id\times r)^*D_A(\varphi)^*\rho_{A*}\hat\omega_{\Gamma}\Bigr]\\
		&= (-1)^{a+1} \sum_{\Gamma} [\Gamma ]\otimes (\id\times r)^*D_A(\varphi)^*\rho_{A*}\hat\omega_{\Gamma}
	\end{split}
\end{equation}
where we have used in the third equality the fact that the form
\begin{equation}\label{eq:condition_alpha}
 \sum_{\Gamma \atop \text{labelled}} [\Gamma] \otimes {\rho_A}_*\hat\omega_{\Gamma} \in
 \calA_k \otimes \Omega^{(n-j-2)k+j+1}_{DR}(\calI_j (\R^n ))
\end{equation}
is closed (the proof of this fact is exactly the same as \cite[Lemma~6.5.15]{R}). Moreover the vanishing of the infinite face contribution together with the generalized Stokes theorem implies that $d\int_{C_1(\R^j)}\beta_\Gamma=\int_{C_1(\R^j)}d\beta_\Gamma$. 
\end{proof}
This completes the proof of Theorem~\ref{thm:closed}(2).

\begin{proof}[Proof of Theorem~\ref{thm:closed2}]
$\calI_j (\R^n)$ is homotopy equivalent to the Stiefel manifold $V_j (\R^n)$ (a deformation retraction is given by the Gram--Schmidt orthogonalization, see e.g., \cite[\S 2.5]{R}) and $\dim V_j(\R^n)=j(2n-j-1)/2$.
Thus, if $k$ is large enough as required in Theorem~\ref{thm:closed2}, then $H^{(n-j-2)k+j+1}(\calI_j(\R^n);\R)=0$. Hence there exists a form $\alpha_k \in \calA_k \otimes \Omega^{(n-j-2)k+j}_{DR}(\calI_j (\R^n ))$ such that $d\alpha_k$ is equal to
(\ref{eq:condition_alpha}).
Then by the definition of $\beta_\Gamma$ and by the generalized Stokes theorem (\ref{eq:stokes}), the correction term is equal to
\[\begin{split}
	\sum_\Gamma[\Gamma]\otimes\int_{C_1(\R^j)}\beta_\Gamma
	&=-\int_{C_1(\R^j)}\mathrm{pr}_{23*}d\widetilde{D}_A^*\alpha_k\\
	&=-d\int_{C_1(\R^j)}\mathrm{pr}_{23*}\widetilde{D}_A^*\alpha_k-Jr^*\int_{C_1(\R^j)}{D}_A(\varphi)^*\alpha_k.
\end{split}
\]
where we have used the fact that $\widetilde{D}_A$ is the constant map near $([0,1]\times \partial C_1(\R^j))\cup (\{0\}\times C_1(\R^j))$. By putting $\bar\alpha_k=(-1)^{j+1}J D_A(\varphi)^*\alpha_k$, we get the result.
\end{proof}

As a consequence of Theorem~\ref{thm:nontrivial} and Proposition~\ref{prop:calA}, $[\hat{z}_k]$ will give a
nontrivial cohomology class for odd $k \ge 3$.
If $k$ is odd and large enough, then $[\bar{z}_k ]$ is also a nontrivial cohomology class of $\emb{n}{j}$ since
$ [\hat{z}_k ] = r^*[\bar{z}_k ]$.

\begin{Rem}\label{rem:anomaly-SH}
It is known that the image of the natural map
\[ f:\pi_0(\emb{5}{3})\to \pi_0(\Imm(\R^3,\R^5))\cong \Z \]
(the isomorphism on the right is given by Smale's isomorphism \cite{Sm}) is $24\Z$. (See \cite{Ek, HM} etc.) We denote this map by $SH:\pi_0(\emb{5}{3})\to 24\Z$. The target of $SH$ is the set of regular homotopy classes of embeddings. It follows from \cite[Theorem~2.5]{B} that the map $SH$ agrees with the composition
\[ \xymatrix{
	\pi_0(\emb{5}{3}) \ar[r]^-{B} & \pi_0(\Omega\Imm(\R^2,\R^4)) \ar[r]^-{G} & \pi_0(\Imm(\R^3,\R^5))\cong \Z 
}\]
of some two maps defined in \cite[Theorem~2.5, Proposition~3.2]{B}.

On the other hand the anomaly correction term $\Theta_k$ defined above gives a 0-form on $\femb{4}{2}$. It is easy to see that when both $n$ and $j$ are even the pullback $i^*\Theta_k$ of $\Theta_k$ by the natural map $i:\Omega\Imm(\R^j,\R^n)\to \femb{n}{j}$ is closed on $\Omega\Imm(\R^j,\R^n)$ and hence gives a well-defined homomorphism
\[ A_k=i^*\Theta_k:\pi_0(\Omega\Imm(\R^2,\R^4))\to \R. \]
At present we do not know the answer to the following question.
\begin{Ques}
Can the map $A_k\circ B:\pi_0(\emb{5}{3})\to \R$ recover $SH$? In other words, is there a non-zero real constant $\lambda_k$ such that $SH=\lambda_k\cdot A_k\circ B$?
\end{Ques}
\end{Rem}

\def\blue#1{{\color{blue}#1}}
\section{Non-triviality of $\hat{z}_k$}\label{sec_ribbon}

Here we will construct the `wheel-like' cycles and evaluate the cohomology classes
$[z_k]\in H^{k(n-j-2)}_{DR}(\emb{n}{j};\calA_k)$ or $[\hat{z}_k]\in H^{k(n-j-2)}_{DR}(\femb{n}{j};\calA_k)$ on the cycles
to show that they are nontrivial for some $k$.

\subsection{Long embeddings from wheel-like ribbon presentations and their special family}

\begin{figure}[htb]
\includegraphics[scale=0.6]{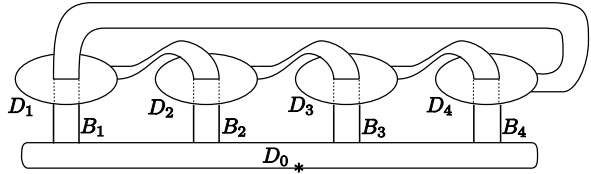}
\caption{The wheel-like ribbon presentation of order $k=4$}\label{wheel_like_ribbon}
\end{figure}

\begin{figure}[htb]
\includegraphics[scale=0.55]{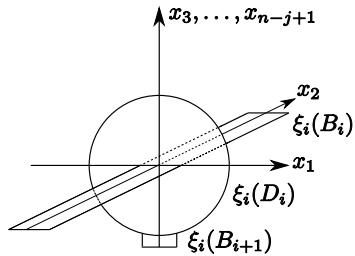}
\caption{Local model of an intersection}\label{local_model}
\end{figure}

\begin{Def}\label{def:wheel-like-pres}
A {\em wheel-like ribbon presentation} $P=D\cup B$ of order $k$ is a based, oriented, immersed $2$-disk in $\R^{n-j+1}$
as shown in Figure~\ref{wheel_like_ribbon}.
More precisely, $P$ consists of $k+1$ disjoint 2-disks
$D=D_0\cup D_1\cup\dots\cup D_k$
and of $k$ disjoint bands
$B=B_1\cup B_2\cup\dots\cup B_k$
($B_i \approx I \times I$ for each $i$), such that
\begin{itemize}
\item
 $B_{i+1}$ connects $D_0$ with $D_i$ ($1\le i\le k$, where $B_{k+1} := B_1$) so that
 $B_{i+1} \cap D_0 = \{ 0\} \times I$, $B_{i+1} \cap D_i = \{ 1\} \times I$,
\item
 each disk $D_i$ intersects `quasi-transversally' with the band $B_i$, $1 \le i \le k$, that is, 
 the intersection $D_i\cap B_i$ is a segment contained in $\Int D_i$ and $TD_i+TB_i$ spans a 3-dimensional subspace at each point in $D_i \cap B_i$ (as in Figure~\ref{local_model}),
\item the base point $*$ of $P$ is on the boundary of $D_0$ but not on the boundaries of $B_i$'s.\qed
\end{itemize}
\end{Def}
Figure~\ref{local_model} shows an image of a neighborhood $U_i$ of $D_i$ via a local homeomorphism $\xi_i : U_i \xrightarrow{\approx} [-3,3]^{n-j+1}$.


\begin{Def}\label{def:V_P}
Define a {\em ribbon $(j+1)$-disk} $V_P$ by
\begin{equation}\label{V_P}
 V_P := \Bigl( D\times\Bigl[ -\frac{1}{2},\frac{1}{2} \Bigr]^{j-1} \Bigr) \cup
 \Bigl( B \times \Bigl[-\frac{1}{4},\frac{1}{4} \Bigr]^{j-1} \Bigr) \subset \R^{n-j+1}\times \R^{j-1}.
\end{equation}
$V_P$ is an immersed handlebody obtained by attaching $1$-handles to $0$-handles in such a way as $P$ indicates, so we can make $V_P$ an immersed $(j+1)$-manifold without corners in the standard way (see e.g.~\cite{K}).
The boundary of $V_P$ is a smoothly embedded $j$-sphere.
Taking a connect-sum of $\partial V_P$ with standard $j$-plane $\iota (\R^j)\subset\R^n$ at the base point, we obtain an embedded $j$-plane in $\R^n$ which is standard
outside a $j$-disk.
We choose a parametrization $\R^j\to\iota (\R^j)\sharp\partial V_P$ for the $j$-plane to obtain a long embedding $\varphi_k:\R^j\hookrightarrow\R^n$.\qed
\end{Def}

\subsubsection{`Resolved' cycles $c_k$, $\widetilde{c}_k$}\label{subsubsec_unclasping_cycle}

\begin{figure}[tb]
\includegraphics[scale=0.55]{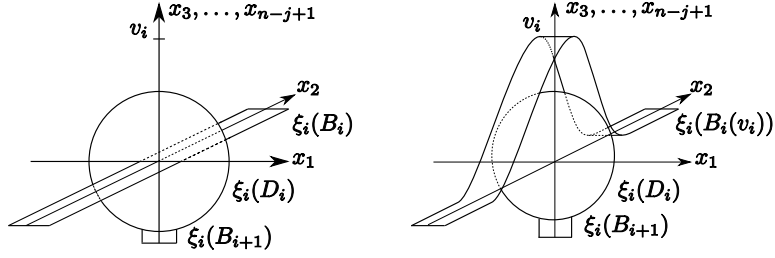}
\caption{Perturbation of a crossing}\label{local_perturbation}
\end{figure}

Here we construct a cycle $c_k$ of $\emb{n}{j}$ of degree $k(n-j-2)$ by `perturbing' the long embedding $\varphi_k$
around the crossings of $\varphi_k$ (neighborhoods of $D_i$'s).
This cycle is a generalization of a `$k$-scheme' in \cite{HKS, Wa},

Consider an $(n-j-2)$-dimensional unit sphere in $x_3 \dots x_{n-j+1}$-space
\[
 S := \{ (0,0,x_3 ,\dots ,x_{n-j+1}) \, \Bigl| \, (x_3 -1)^2 + x_4^2 +\dots +x^2_{n-j+1}= 1 \} .
\]
We perturb $B_i$ by considering, for any $v \in S$, a (2-dimensional) band
\[
 B(v):= \Bigl\{ (x,y;\gamma (y) v) \in \R^2 \times \R^{n-j-1} \, \Bigl| \, \abs{x} \le \frac{1}{2},\ \abs{y} < 3
 \Bigr\}
\]
(see Figure \ref{local_perturbation}) where $ \gamma (y) := \exp \bigl(- {y^2}/{\sqrt{9-y^2}} \bigr)$.
Replacing each $B_i$ by $B_i (v_i ) :=(B_i \setminus (B_i \cap U_i )) \cup \xi^{-1}_i (B(v_i ))$, we obtain a new
ribbon presentation $P_{\mathbf{v}} := D \cup B_{\mathbf{v}}$ for any $\mathbf{v} := (v_1 ,\dots ,v_k )\in (S^{n-j-2})^{\times k}$, where $B_{\mathbf{v}} := B_1 (v_1 )\cup\dots\cup B_k (v_k )$.
Taking the boundary of the $(j+1)$-disk $V_{P_{\mathbf{v}}}$, we have a long embedding $\varphi_k^{\bf v}$,
a `perturbation' of $\varphi_k$ via $\mathbf{v} \in (S^{n-j-2})^{\times k}$.
We can take $\varphi^{\bf v}_k$ to be continuous with respect to ${\bf v}$ (see the remark below).
Thus we have a continuous map
\[
 c_k : (S^{n-j-2})^{\times k} \longrightarrow \emb{n}{j},\quad {\bf v} \longmapsto \varphi_k^{{\bf v}}.
\]
This is canonical up to homotopy.
We regard the map as a $k(n-j-2)$-cycle of $\emb{n}{j}$.

Moreover, we have not only a family of embeddings but also a family $\{V_{P_{{\bf v}}}\}_{\bf v}$ of ribbon disks.
We get a family of paths in $\Imm(\R^j,\R^n)$
\[
 [0,1]\times (S^{n-j-2})^{\times k}\longrightarrow \Imm (\R^j ,\R^n )
\]
such that each path in this family
collapses
each embedding $\varphi^{\bf v}_k$ (${\bf v}\in (S^{n-j-2})^{\times k}$) to the standard inclusion along the ribbon disk $V_{P_{\bf v}}$ by a regular homotopy.
Inverting each path, we obtain
a map $\widetilde{c}_k:(S^{n-j-2})^{\times k}\to \femb{n}{j}$ which extends $c_k$.
We will consider $\widetilde{c}_k$ as representing a cycle of $\femb{n}{j}$.

\begin{Rem}\label{rem:take_emb}
A reason why it is possible to take a family of embeddings $c_k$ for the family of submanifolds $\{\partial V_{P_{\mathbf{v}}}\}_{\mathbf{v}}$ is that the relative smooth $(\R^j,\R^j\setminus D^{j})$-bundle over $(S^{n-j-2})^{\times k}$ given by the family $\{\partial V_{P_{\mathbf{v}}}\}_{\mathbf{v}}$ is trivial because it can be collapsed to a constant family that is isotopic to the standard inclusion by a sequence of unclaspings on every crossings that are given through a family of immersions. 

The support of the deformation can be restricted inside the union of the crossings. Thus we may assume that the family $\{\varphi_k^{\mathbf{v}}\}_{\mathbf{v}}$ is constant outside crossings.\qed
\end{Rem}

\subsubsection{Main evaluation}

\begin{figure}[htb]
\includegraphics{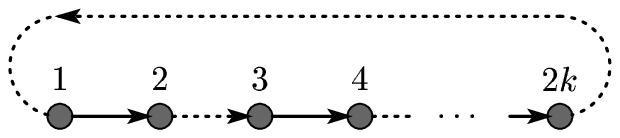}
\caption{The polygonal graph $\Gamma^{(k)}$}\label{graph_with_no_hair}
\end{figure}

Let $\Gamma^{(k)}$ be the polygonal graph defined by Figure~\ref{graph_with_no_hair}.
In the rest of this section, we will prove the following theorem.

\begin{Thm}\label{thm:nontrivial}
\begin{enumerate}
\item Suppose $n,j,k$ are as in Theorem~\ref{thm:closed} (1); (a) $n$ odd, or (b) $n$ even, $j$ odd and $k\le 4$, or (c) $n\ge 12$ even, $j=3$.
Then $\pair{z_k}{c_k}=\pm [\Gamma^{(k)}]$, where $\pair{\alpha}{c_k}$ denotes $\int_{(S^{n-j-2})^{\times k}}c_k^*\alpha$.
 Thus both $[{z}_k]\in H^{k(n-j-2)}_{DR}(\emb{n}{j};\calA_k)$ and $[{c}_k]\in H_{k(n-j-2)}(\emb{n}{j};\R)$ are nontrivial
if $k\ge 2$ is such that $[\Gamma^{(k)}]\neq 0$ in $\calA_k=\calA_k(n,j)$.
\item If $n,j$ are both even as in Theorem~\ref{thm:closed} (2), then $\pair{\hat{z}_k}{\widetilde{c}_k}=\pm [\Gamma^{(k)}]$.
 Thus both $[\hat{z}_k]\in H^{k(n-j-2)}_{DR}(\femb{n}{j};\calA_k)$ and $[\widetilde{c}_k]\in H_{k(n-j-2)}(\femb{n}{j};\R)$ are nontrivial if $[\Gamma^{(k)}]\ne 0$ in $\calA_k$.
If moreover $n\ge 2j$, then
$r_*[\widetilde{c}_k]\in H_{k(n-j-2)}(\emb{n}{j};\R)$ is also nontrivial, where $r:\femb{n}{j} \to \emb{n}{j}$ is the forgetting map.
\end{enumerate}
\end{Thm}

\begin{Rem}\label{rem:new_class}
What we know about the space $\calA_k$ are summarized in Proposition~\ref{prop:calA} which will be proved in \S\ref{sec:calA}.
In particular we will show that $[\Gamma^{(3)}]\ne 0$ in
$\calA_3 \cong \R$ when $n$ is odd and $j$ is even (Proposition~\ref{prop:A3_odd_codim}).
Hence by Theorem~\ref{thm:nontrivial} (1), $[z_3 ] \in H^{3(n-j-2)}_{DR}(\emb{n}{j})$ is not zero.
To the authors' knowledge, this is the first cohomology class of higher degree than the homology classes discussed in \cite{B} (in the cases where $n$ is odd and $j$ is even).\qed
\end{Rem}

The proof is outlined as follows.
We may compute $\pair{z_k}{c_k}$ or $\pair{\hat{z}_k}{\widetilde{c}_k}$ in the limit that the crossings of $\varphi_k$
`shrink to a point' (see \S\ref{sss:shrinking}) since a shrinking of a crossing does not change $[c_k], [\widetilde{c}_k]$ and since $z_k,\hat{z}_k$ are closed.
We will show in \S\ref{ss:evaluation} that, in the limit,
\[
 \pair{I(\Gamma )}{c_k} \longrightarrow
 \begin{cases}
  \pm \abs{\mathrm{Aut}\, \Gamma} & \text{if } \Gamma =\Gamma^{(k)} \text{ polygonal with no} \\
	& \text{ orientation reversing automorphism},\\
  0 & \text{otherwise}
 \end{cases}
\]
and that the value of the correction term for $\hat{z}_k$ on $\widetilde{c}_k$ vanishes when $n, j$ are even.
Here $\Aut\Gamma$ denotes the automorphism group of the underlying (unoriented) graph $\Gamma$.
Since the polygonal graph is unique for each $k$, the pairing $\pair{z_k}{c_k}$($=\pair{\hat{z}_k}{\widetilde{c}_k}$ when
$n,j$ are even) is equal to $\pm [\Gamma^{(k)}]$.

\subsection{Modification of embeddings to convenient ones}

For the convenience in evaluating the integral, we deform the family $c_k=\{\varphi_k^\mathbf{v}\}_\mathbf{v}$ (keeping the property mentioned in Remark~\ref{rem:take_emb} satisfied) as follows.

\subsubsection{Shrinking}\label{sss:shrinking}
\begin{figure}[htb]
\includegraphics[scale=0.55]{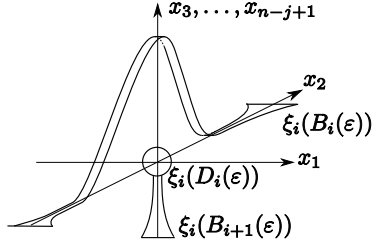}
\caption{A shrinking of the crossing (compare it with Figure~\ref{local_perturbation})}\label{shrink_crossing}
\end{figure}

Let $\varepsilon >0$ be sufficiently small.
We choose a ribbon presentation $P$ so that the neighborhoods $U_i =\xi^{-1}_i ([-3,3]^{n-j+1})$ of the crossings of
$\varphi_k$ are contained in $\varepsilon$-balls.
We also deform the local model of the crossings of $\varphi^{\mathbf{v}}_k$ as in Figure~\ref{shrink_crossing}, replacing the bands and the disks with
\begin{align*}
 B(\varepsilon ) &:= 
 \Bigl\{ (x,0,z,\mathbf{0}) \in \R^3 \times \{ 0\}^{n-j-2} \, \Bigl| \,
  -3 \le z \le -\sqrt{\varepsilon^2 -x^2},\ \abs{x} \le \frac{\varepsilon^2}{2}
  \Bigr\} , \\
 D(\varepsilon ) &:= \{ (x,0,z,\mathbf{0}) \in \R^3 \times \{ 0\}^{n-j-2} \, | \, x^2 +z^2 \le \varepsilon^2 \}
\end{align*}
and for any $v \in S$,
\[
 B (v ,\varepsilon ) := \Bigl\{ (x,y, \gamma (y)v) \, \Bigl| \,
 \abs{x} \le \frac{1}{2} - \frac{1-\varepsilon^2}{2}\gamma (y),\ \abs{y} < 3 \Bigr\}
\]
(recall $\gamma (y) = e^{-y^2 / \sqrt{9-y^2}}$).
Replacing $D_i \cap U_i$, $B_{i+1} \cap U_i$ and $B_i (v_i ) \cap U_i$ with
\[
 D_i (\varepsilon ) := \xi^{-1}_i (D(\varepsilon )),
 \qquad B_{i+1} (\varepsilon ) := \xi^{-1}_i (B(\varepsilon )),
 \qquad B_i (v_i ,\varepsilon ) := \xi^{-1}_i (B(v_i ,\varepsilon )),
\]
we obtain a new perturbation of the ribbon presentation, which we denote by
$P_{\mathbf{v},\varepsilon}:= D_{\mathbf{v},\varepsilon} \cup B_{\varepsilon}$.
Then we `fatten' $P_{\mathbf{v},\varepsilon}$ in a similar way to (\ref{V_P}) to obtain $V_{P_{\mathbf{v},\varepsilon}}$,
but now around $U_i$ we fatten $D_i (\varepsilon )$ and $B_i (v_i ,\varepsilon )$ by
$[-\varepsilon /2,\varepsilon /2]^{j-1}$ and $[-\varepsilon^2 /4,\varepsilon^2 /4]^{j-1}$ respectively.
Taking the boundary of $V_{P_{\mathbf{v},\varepsilon}}$, we obtain a family of long embeddings denoted by
$\varphi^{\mathbf{v}, \varepsilon}_k$.

Clearly the choice of $\varepsilon \in (0,1)$ does not affect the homology classes $[c_k],[\widetilde{c}_k]$.
So it is enough to compute $\pair{z_k}{c_k}$ in the limit $\varepsilon \to 0$.

\subsubsection{Crossing as embeddings from standard disks}
\begin{Def}[Crossing]
We write $\hat{U}_i:= U_i\times [-3/4,3/4]^{j-1}$.
Then the intersection of $\hat{U}_i$ with the image of the long embedding $\varphi^{\mathbf{v}, \varepsilon}_k$ separates
into two components.
We denote them by $\hat{D}_i(\varepsilon )\cup\hat{B}_i(v_i,\varepsilon )$, where the two components correspond respectively to $D_i$ and $B_i$.
We call the triple $(\hat{U}_i ,\hat{D}_i (\varepsilon ), \hat{B}_i (v_i ,\varepsilon ))$ the {\it $i$-th crossing} of
$\varphi^{\mathbf{v}, \varepsilon}_k$.\qed
\end{Def}

\begin{figure}[htb]
\includegraphics{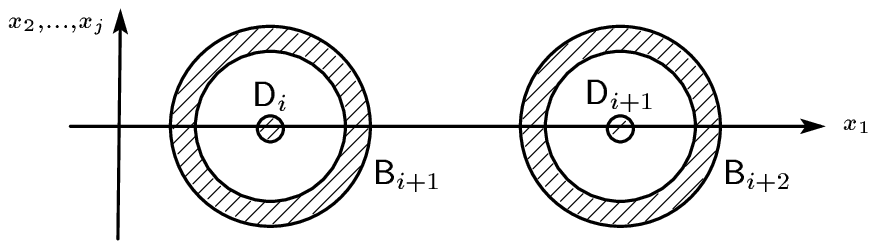}
\caption{$\mathsf{D}_i$ and $\mathsf{B}_{i+1}$}\label{domains}
\end{figure}

$\hat{D}_i (\varepsilon )$ is diffeomorphic to a punctured $j$-sphere and $\hat{B}_i (v_i ,\varepsilon )$ is
diffeomorphic to $I \times S^{j-1}$.
After a suitable deformation, we may assume that, for any $\mathbf{v} \in (S^{n-j-2})^{\times k}$, the parametrization
$\varphi^{\mathbf{v}, \varepsilon}_k :\R^j \hookrightarrow \R^n$ is chosen so that
$\mathsf{D}_i =\mathsf{D}_i (\varepsilon )$, $\mathsf{B}_i =\mathsf{B}_i (\varepsilon )$ are mapped homeomorphically onto
$\hat{D}_i (\varepsilon )$ and $\hat{B}_i (v_i ,\varepsilon)$ respectively, where
\[\begin{split}
	\mathsf{D}_i (\varepsilon ) &:= \{ (x_1 ,\dots ,x_j ) \in \R^j \, | \,
 ( x_1 -p_i) ^2 + x^2_2 + \dots +x^2_j \le (\varepsilon^2 )^2 \}\\
 	\mathsf{B}_i (\varepsilon ) &:= \{ (x_1 ,\dots ,x_j ) \in \R^j \, | \,
  (3\varepsilon/4)^2 \le ( x_1 -p_{i-1}) ^2 + x^2_2 + \dots +x^2_j  \le \varepsilon^2 \} ,
\end{split}
\]
and where $p_i=i/k$ ($1\le i\le k-1$), $p_k=p_0=0$  (see Figure~\ref{domains}).

\subsection{Evaluation by $z_k$}\label{ss:evaluation}
Here we give a proof of Theorem \ref{thm:nontrivial}.
We work with the assumptions on $c_k=\{\varphi_k^{\mathbf{v},\varepsilon}\}_\mathbf{v}$ made in the previous subsection.

\subsubsection{Non-corrected case; $n,j,k$ are as in Theorem~\ref{thm:closed} (i)}\label{ss:non_correction}
From now on we compute the value of
\[
 \pair{z_k}{c_k}
 = \frac{1}{k_S!k_T!}\sum_{{\Gamma}\atop{\mathrm{labelled}}}[\Gamma ] \pair{I(\Gamma )}{c_k}
 = \sum_{{\Gamma}\atop{\mathrm{unlabelled}}} \frac{[\Gamma]}{\abs{\mathrm{Aut}\, \Gamma}}\pair{I(\Gamma )}{c_k},
\]
where in the last term $\Gamma$ runs over all unlabelled admissible 1-loop graphs of order $k$ and where $I(\Gamma)$ and $[\Gamma]$ are given for some common labelled representative for each unlabelled graph $\Gamma$.
Note that  there are $\frac{k_S!k_T!}{|\Aut \Gamma|}$ different labellings on a graph $\Gamma$ and that the product $[\Gamma]I(\Gamma)$ does not depend on the choice of a label.
We compute each term $\pair{I(\Gamma )}{c_k}$ explicitly for all $\Gamma$.

Let $s = \abs{V_i (\Gamma )}$, $t=\abs{V_e (\Gamma )}$.
Consider the following commutative diagram;
\[
 \xymatrix{ 
  (S^{n-j-2})^{\times k} \times C_s (\R^j ) \ar[r]^-{c_k \times \mathrm{id}} \ar[d]_-{\mathrm{pr}_1}
   & \emb{n}{j} \times C_s (\R^j ) \ar[d]^-{\mathrm{pr}_1} & \ar[l]_-{p_\Gamma} \ar[ld]^-{\pi_{\Gamma}}C_\Gamma \\
  (S^{n-j-2})^{\times k} \ar[r]^-{c_k} & \emb{n}{j}
 }
\]
where $p_{\Gamma}$ is given by $(\varphi;x_1,\ldots,x_s;x_{s+1},\ldots,x_{s+t})\mapsto (\varphi;x_1,\ldots,x_s)$.
Then
\[
 \pair{I(\Gamma )}{c_k} = \int_{(S^{n-j-2})^{\times k}}c^*_k (\pi_{\Gamma})_* \omega_{\Gamma} \\
 = \int_{(S^{n-j-2})^{\times k} \times C_s (\R^j )} (c_k \times \mathrm{id})^* (p_{\Gamma})_* \omega_{\Gamma}.
\]

\begin{Lem}\label{V1_contribution}
Let $V_1(i)$ be the subset of $C_s (\R^j )$ consisting of configurations such that at most one point of a configuration is in $\mathsf{D}_i (\varepsilon )\cup \mathsf{B}_i (\varepsilon )$.
Then
\[ \int_{(S^{n-j-2})^{\times k} \times V_1(i)}
 (c_k \times \mathrm{id})^* (p_{\Gamma})_* \omega_{\Gamma} =O(\varepsilon )
\]
(this means that the left hand side converges to zero as $\varepsilon$ tends to zero).
\end{Lem}

\begin{proof}
If one of $\mathsf{D}_i$ and $\mathsf{B}_i$ contains no points, then the integral differs only by $O(\varepsilon )$ from an integral of a pullback of a $k(n-j-2)$-form on $(S^{n-j-2})^{\times k-1}\subset (S^{n-j-2})^{\times k}$ (the complemental direction of the $i$-th factor) along the projection.
This is because we can deform $c_k$ in $\hat{U}_i$, by a small regular homotopy, so that $c_k$ is constant for any $v_i \in S^{n-j-2}$ and the integral remains to be well-defined all through the deformation.
The integral changes only by $O(\varepsilon )$ since the change of $\phi_e$ (regarded as a smooth map from $C_{\Gamma}\times (S^{n-j-2})^{\times k}$) by the deformation can be made arbitrarily small.
\end{proof}

The pairing $\pair{z_k}{c_k}$ is independent of the choice of $\varepsilon$ since the homology class $[c_k]$ is independent of $\varepsilon$ and $z_k$ is closed by the assumption on $n,j$.
Thus by Lemma \ref{V1_contribution} we may restrict to the integration on the subspace of $C_s (\R^j )$
consisting of configurations such that at least one point is mapped to both $\hat{D}_i$ and $\hat{B}_i$ by $\varphi_k$
(other configurations contribute to the integral by $O(\varepsilon )$).

Since $c_k$ has exactly $k$ crossings $(\hat{U}_i ,\hat{D}_i , \hat{B}_i )$, $\Gamma$ has to satisfy $s \ge 2k$ to
contribute to the pairing $\pair{z_k}{c_k}$ nontrivially in the limit $\varepsilon \to 0$.
But since $\Gamma$ is of order $k$, we have $s+t=2k$ vertices (Definition~\ref{def:order}) and thus $s \le 2k$.
Hence only the graphs with $s=2k$ (and thus $t=0$, that is, without e-vertices) can contribute nontrivially to the
pairing $\pair{z_k}{c_k}$.

\begin{Lem}\label{V2_contribution}
Let $\Gamma$ be an admissible graph without e-vertices, and $e=\overrightarrow{pq}$ its $\eta$-edge.
Let $V_2 (e)$ be the subspace of $C_{\Gamma}(\varphi )\cong C_{2k}(\R^j )$ consisting of configurations such that
the points corresponding to $p$ and $q$ are not in the same $\mathsf{S}_i$, where $\mathsf{S}_i$ is a $j$-ball containing
$\mathsf{D}_i \cup \mathsf{B}_{i+1}$ $(\mathsf{B}_{k+1}:=\mathsf{B}_1)$ ;
\[
 \mathsf{S}_i := \{ (x_1 ,\dots ,x_j ) \in \R^j \, | \, ( x_1 -p_i) ^2 + x^2_2 +\dots +x^2_j \le \varepsilon^2 \}
\]
where $p_i=i/k$ $(1\le i\le k-1)$, $p_k=0$.
Then
\[
 \int_{(S^{n-j-2})^{\times k} \times V_2 (e)} (c_k \times \mathrm{id})^* (p_{\Gamma})_* \omega_{\Gamma}
 = O(\varepsilon ).
\]
\end{Lem}

\begin{proof}
By Lemma \ref{V1_contribution}, only the configurations where each one of $2k$ points belongs to one $\mathsf{S}_i$ can
contribute nontrivially to $\pair{z_k}{c_k}$.
If the points $x_p$ and $x_q$ are in different $\mathsf{S}_i$'s, then the image of the map $\phi_e$ concentrates in some small ball (with radius $O(\varepsilon )$) in $S^{j-1}$, because of the assumption for
$\mathsf{D}_i (\epsilon )$ and $\mathsf{B}_i (\epsilon )$.
Thus the integral of a product of edge forms over $V_2 (e)$ is $O(\varepsilon )$.
\end{proof}

\begin{Lem}\label{V3_contribution}
Let $\Gamma$ be an admissible graph without e-vertices, and $e=\overrightarrow{pq}$ its $\theta$-edge.
Let $V_3 (e)$ be the subspace of $C_{\Gamma}(\varphi )\cong C_{2k}(\R^j )$ consisting of configurations with
$(x_p ,x_q ) \not\in {\sf D}_i \times {\sf B}_i$ and $\not\in {\sf B}_i \times {\sf D}_i$ for any $i$.
Then
\[
 \int_{(S^{n-j-2})^{\times k} \times V_3 (e)} (c_k \times \mathrm{id})^* (p_{\Gamma})_* \omega_{\Gamma}
 = O(\varepsilon ).
\]
\end{Lem}

\begin{proof}
By assumption and Lemma~\ref{V1_contribution}, we may assume $(x_p,x_q)\in {\sf D}_i\times{\sf B}_{i'}$ or
$\in {\sf B}_i\times{\sf D}_{i'}$ for some $i\ne i'$.
But then the image of $\phi_e$ is in a small $(n-1)$-disk (of radius $O(\varepsilon )$) in $S^{n-1}$.
\end{proof}

\begin{Lem}\label{wheel-like_contribution}
In the limit $\varepsilon \to 0$,
\[
 \pair{z_k}{c_k} = \pm \frac{[\Gamma^{(k)}]}{\abs{\mathrm{Aut}\, \Gamma^{(k)}}}\pair{I(\Gamma^{(k)})}{c_k}
 + O(\varepsilon ),
\]
where $\Gamma^{(k)}$ is the unique polygonal graph (see Figure \ref{graph_with_no_hair}) of order $k$.
\end{Lem}

\begin{proof}
Let $\Gamma$ be a graph without e-vertices.
If an i-vertex $p$ is trivalent (thus $\Gamma$ is not polygonal), there are two $\eta$-edges (say $pq$ and $pr$)
and one $\theta$-edge emanating from $p$.
Then by the above Lemma \ref{V2_contribution}, the three points $x_p$, $x_q$ and $x_r$ must be in the same
$\mathsf{S}_i$.
But then there must be one $\mathsf{D}_l$ or $\mathsf{B}_l$ which contains no points in a configuration.
Thus for any $\Gamma$ which is not polygonal, we have $\pair{I(\Gamma )}{c_k}=O(\varepsilon )$ by Lemma~\ref{V3_contribution} and by the identity
\[ \bigcup_{e\in E_{\eta}(\Gamma ),e'\in E_{\theta}(\Gamma )} (V_2 (e) \cup V_3 (e'))=C_{2k}(\R^j )\]
for such a graph $\Gamma$.
\end{proof}

The final task is to compute $\pair{I(\Gamma^{(k)})}{c_k}$, where $\Gamma^{(k)}$ is the polygonal graph oriented as in
Figure \ref{graph_with_no_hair}.
We prove the following lemma.
\begin{Lem}\label{lem:I(G)=Aut}
If $k$ is such that the polygonal graph $\Gamma^{(k)}$ does not have an orientation reversing automorphism, then
\[
 \pair{I(\Gamma^{(k)})}{c_k}=\pm \abs{\Aut\Gamma^{(k)}}.
\]
Otherwise $\pair{I(\Gamma^{(k)})}{c_k}=0$.
\end{Lem}

\begin{proof}
By Lemma \ref{V1_contribution}, we may restrict the integration on the configurations where all the points are in one of
$\mathsf{D}$'s or $\mathsf{B}$'s.
By Lemma~\ref{V2_contribution} it suffices to consider only the case where the points $x_{2i-1}$, $x_{2i}$ corresponding to endpoints $2i-1$, $2i$ of an $\eta$-edge must be in $\mathsf{D}_l$ and $\mathsf{B}_{l+1}$ for some $l$.
Then by Lemma~\ref{V3_contribution}, $x_{2i}$ must be in $\mathsf{B}_{l+1}$ (hence $x_{2i-1}\in \mathsf{D}_l$) and
the endpoint $x_{2i+1}$ of a $\theta$-edge other than $x_{2i}$ is forced to be in $\mathsf{D}_{l+1}$.
There are $\abs{\mathrm{Aut}\, \Gamma^{(k)}}=2k$ components of such configurations as above
(because $\Aut\Gamma^{(k)}$ is isomorphic to the dihedral group of the $k$-gon).
By symmetry it is enough to compute the integral on the component $\Pi_k$ of
$C_{2k}(\R^j)\setminus \bigcup_{e,e'}V_2(e)\cup V_3(e')$ among the $2k$ components where the configuration satisfies
$x_{2i-1} \in \mathsf{D}_i$, $x_{2i} \in \mathsf{B}_{i+1}$ $(1 \le i \le k)$.
Other components contribute to the integral by the same value modulo signs as the component $\Pi_k$.
The sign which is induced by a permutation of vertices is the same as that induced on the graph by the corresponding
permutation.
Therefore the integral $\pair{I(\Gamma^{(k)})}{c_k}$ vanishes by self-cancelling if $\Gamma^{(k)}$ has an orientation
reversing automorphism.

We claim that, when $\Gamma^{(k)}$ does  not have an orientation reversing automorphism, the integral $\pair{I(\Gamma^{(k)})}{c_k}$ restricted to $\Pi_k$ is the product of the `linking numbers' of
$\hat{D}_i (\varepsilon )$ with $\bigcup_{v_i \in S} \hat{B}_i (v_i ,\varepsilon )$ ($1 \le i \le k$), which are equal
to $\pm 1$.
We will see this more rigorously now:

To describe $\pair{I(\Gamma^{(k)})}{c_k}$ explicitly, we define two types of direction maps;
\begin{gather*}
 \phi_{\theta ,i} : \mathsf{D}_i \times \mathsf{B}_i \times S^{n-j-2} \longrightarrow S^{n-1}, \quad
 (d_i ,b_i ,v_i )\mapsto u(\varphi^{v_i}_k (d_i )-\varphi^{v_i}_k (b_i )),\\
 \phi_{\eta ,i} : \mathsf{D}_i \times \mathsf{B}_{i+1} \longrightarrow S^{j-1}, \quad
 (d_i ,b_{i+1})\mapsto u(b_{i+1} - d_i),
\end{gather*}
where $d_i \in \mathsf{D}_i$, $b_i \in \mathsf{B}_i$, $\varphi^{v_i}_k$ is the embedding $\varphi_k$ with its $i$-th
crossing perturbed by $v_i$, and $u(v):=v/\abs{v}$ for a nonzero vector $v$.
Then by Lemmas \ref{V1_contribution}, \ref{V2_contribution} and \ref{V3_contribution}, we have
\begin{equation}\label{eq:I(Gk)_ck}
 \pair{I(\Gamma^{(k)})}{c_k}=2k\int_{\Pi_k \times (S^{n-j-2})^{\times k}}
 \bigwedge_{i=1}^k \phi_{\theta ,i}^* vol_{S^{n-1}} \wedge \phi_{\eta ,i}^* vol_{S^{j-1}} + O(\varepsilon ).
\end{equation}
But we can replace $\phi_{\eta ,i}$ (changing the integral (\ref{eq:I(Gk)_ck}) only by $O(\varepsilon )$) by
\[
 \phi^o_{\eta ,i} : \mathsf{B}_{i+1} \longrightarrow S^{j-1}, \quad b_{i+1} \longmapsto u(b_{i+1}),
\]
because our $\mathsf{D}_i$ is quite smaller than $\mathsf{B}_{i+1}$, and consequently \eqref{eq:I(Gk)_ck} can be rewritten as
\begin{align*}
 &\int_{\Pi_k \times (S^{n-j-2})^{\times k}}
  \bigwedge_{i=1}^k \phi_{\theta ,i}^* vol_{S^{n-1}} \wedge (\phi^o_{\eta ,i})^* vol_{S^{j-1}} + O(\varepsilon ) \\
 &\ = \prod_{i=1}^k\int_{\mathsf{D}_i \times \mathsf{B}_i \times S^{n-j-2}}
  \phi_{\theta ,i}^* vol_{S^{n-1}} \wedge (\phi^o_{\eta ,i-1})^* vol_{S^{j-1}} + O(\varepsilon ).
\end{align*}
Then Lemma~\ref{lem:following-lem} below completes the proof of Lemma~\ref{lem:I(G)=Aut}.
\end{proof}

\begin{Lem}\label{lem:following-lem}
\[
 \int_{\mathsf{D}_i \times \mathsf{B}_i \times S^{n-j-2}}
  \phi_{\theta ,i}^* vol_{S^{n-1}} \wedge (\phi^o_{\eta ,i-1})^* vol_{S^{j-1}}=\pm 1+O(\varepsilon).
\]
\end{Lem}

\begin{proof}
Under the identifications $\mathsf{D}_i \approx D^j$ and $\mathsf{B}_i \approx I \times S^{j-1}$,
the map $\phi_{\theta ,i} \times \phi^o_{\eta ,i-1}$ can be seen as
\begin{gather*}
 D^j \times I \times S^{j-1} \times S^{n-j-2}\to S^{n-1}\times S^{j-1}, \qquad
 (x,t,w,v)\mapsto \left( u(\varphi_k (x)-\varphi^v_k (t,w)), w \right)
\end{gather*}
The point $\varphi^v_k (t,u)$ is in the cylinder $\hat{B}_i (v,\varepsilon ) \approx I \times S^{j-1}$, which has as its
`core' an arc
\[
 \gamma (v,t)=\left( 0,t,v\exp ( -t^2/\sqrt{9-t^2})\right)
\]
(see \S \ref{subsubsec_unclasping_cycle}), and is fattened by taking a product with a small $S^{j-1}$ in $x_1 x_{n-j+2} \dots x_n$-direction.
Since the radius of the $S^{j-1}$ is quite smaller ($\sim \varepsilon^2$) than that of $\hat{D}_i (\varepsilon )$ ($\sim \varepsilon$), the map $\phi_{\theta ,i}$ can be replaced (changing the integral only by $O(\varepsilon )$) by the map
\[
 \phi^o_{\theta,i}: \mathsf{D}_i \times \mathsf{B}_i \times S^{n-j-2} \to S^{n-1}, \qquad
 (x,t,v) \mapsto u(\varphi_k (x)- \gamma (v,t)).
\]
Thus the integral of the statement is rewritten as
\[
 \int_{S^{j-1}}(\phi^o_{\eta,i-1})^* vol_{S^{j-1}}
 \int_{D^j \times I \times S^{n-j-2}}(\phi^o_{\theta ,i})^* vol_{S^{n-1}} +O(\varepsilon ).
\]
The first integral is obviously one, since $\phi^o_{\eta ,i-1}$ restricts to the identity on $S^{j-1}$.
The second integral is $lk(A_i ,\mathcal{S})+O(\varepsilon )$, where $lk$ is the linking number,
\[
 A_i := \bigcup_{t \in I}\bigcup_{v_i \in S^{n-j-2}}\gamma (v_i ,t) \approx \Sigma S^{n-j-2},
\]
and $\mathcal{S}$ is a $j$-sphere obtained from $\hat{D}_i (\varepsilon )$ by stopping up a small $j$-ball (corresponding
to $D_i \cap B_{i+1}$).
$\mathcal{S}$ is a unit $j$-sphere in $x_2 x_3 x_{n-j+2}\dots x_n$-space centered at the origin, and
$A_i$ is a unit $(n-j-1)$-sphere in $x_1 x_3x_4 \dots x_{n-j+1}$-space centered at $(0,0,1,0,\dots ,0)$.
Thus $lk(A_i,\mathcal{S})$ is clearly $\pm 1$.
\end{proof}

Lemmas~\ref{wheel-like_contribution}, \ref{lem:I(G)=Aut} complete the proof of Theorem~\ref{thm:nontrivial} (1).

\subsubsection{The correction term\blue{;} $n,j$ are even}\label{ss:correction-change}
In the case where $n,j$ are both even, instead of evaluating $\pair{\hat{z}_k}{\widetilde{c}_k}$, we compute the difference
\[ \pair{\hat{z}_k}{\widetilde{c}_k} - \pair{\hat{z}_k}{\widetilde{c}_k^0} \]
for some nullhomotopic cycle $\widetilde{c}_k^0$ of $\femb{n}{j}$ given as follows. 

\begin{figure}[htb]
\fig{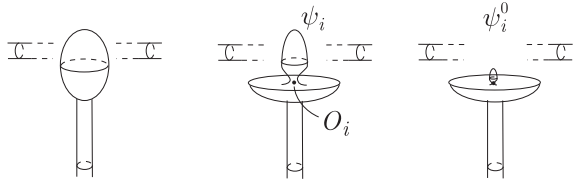}
\caption{Unclasping by scaling down around a point $O_i$.}\label{fig:psi_dilation}
\end{figure}

Let $\psi_i$ denote the restricted embedding $c_k(\mathbf{v}^0)|_{\mathsf{D}_i}:\mathsf{D}_i\hookrightarrow \R^n$ where $\mathbf{v}^0\in (S^{n-j-2})^{\times k}$ is the basepoint. Let $O_i$ be the center of the $j$-disk $\partial \hat{U}_i\cap (D_i\times[-\varepsilon/2,\varepsilon/2]^{j-1})$ and fix a local coordinate around $O_i$ induced from that of $\R^n$ so that $O_i$ is the origin. After a suitable deformation of $c_k(\mathbf{v})|_{\mathsf{D}_i}$, we may assume that $\psi_i$ agrees with the standard linear inclusion $\iota$ on $r_0\le \abs{x}\le r_1$ for some $r_0$, $r_1$ with $r_0/r_1\ll 1$, with respect to the local coordinate. Then we set
\[ \psi_i^0(x)=\left\{
	\begin{array}{ll}
		\lambda \psi_i(\lambda^{-1}x) & |x|\leq r_0\\
		\psi_i(x) &  r_0< |x|\leq r_1
	\end{array}\right.
\]
under the local coordinate, for a small constant $\lambda>0$ such that $r_0/r_1 <\lambda <1$, which implies $r_0 <r_0 /\lambda <r_1$. See Figure~\ref{fig:psi_dilation}.
We may also assume that if $\lambda$ is small enough, then the $(j+1)$-disk $D_i (\varepsilon )\times [-\varepsilon/2,\varepsilon/2]^{j-1}$ (after a suitable deformation) does not intersect $\hat{B}_i(v_i,\varepsilon)$ for all $v_i\in S^{n-j-2}$. The resulting embedding $\psi_i^0$ has the same differential $D\psi_i^0:\mathsf{D}_i\to \calI_{j}(\R^n)$ as $\psi_i$ up to a relative isotopy of the domain $\mathsf{D}_i$. More precisely, by definition the differential of $\psi_i^0$ is
\[ D\psi_i^0(x)=\left\{
	\begin{array}{ll}
		D\psi_i(\lambda^{-1}x) & |x|\leq r_0\\
		D\psi_i(x)\ (=\iota) &  r_0< |x|\leq r_1
	\end{array}\right.
\]
Note that this is continuous because $\psi_i$ is standard on $r_0\le \abs{x}\le r_1$.
We deform $\psi_i^0$ by a relative isotopy of $(\mathsf{D}_i ,\partial \mathsf{D}_i )$ so that $D\psi_i^0$ coincides with $D\psi_i$ (we will denote the resulting embedding again by $\psi_i^0$).
Replacing $\psi_i$ with $\psi_i^0$ for all $i$, we get a family of homotopies through immersions
\[ \widetilde{c}_k^0:(S^{n-j-2})^{\times k}\to \femb{n}{j} \]
with the following properties:

\begin{Lem}\label{lem:ck0}
\begin{enumerate}
\item The correction terms evaluated on $\widetilde{c}_k$ and $\widetilde{c}_k^0$ coincide.
\item $\widetilde{c}_k^0$ is nullhomotopic.
\item $\pair{z_k}{c_k^0}=0$ where $c_k^0=r\circ \widetilde{c}_k^0$.
\end{enumerate}
\end{Lem}

\begin{proof}
(2) is because the family $\widetilde{c}_k^0$ of homotopies is in fact a family of embeddings of $[0,1]\times \R^j$.
(3) is checked by the same argument as in the computation of $\pair{z_k}{c_k}$; $\widetilde{c}^0_k$ is arranged so that
the linking numbers of Lemma \ref{lem:following-lem} are zero.
(1) is proved as follows.
The correction term is defined as in \S\ref{ss:anomaly_correction} and its value on $\widetilde{c}_k$ is given by
\[
 \sum_{\Gamma}[\Gamma ]\otimes \int_{[0,1]\times C_1 (\R^j )\times (S^{n-j-2})^{\times k}}
  \hat{D}^* \rho_{V(\Gamma )*}\hat{\omega}_{\Gamma}
\]
where $\hat{D}:[0,1]\times C_1 (\R^j )\times (S^{n-j-2})^{\times k}\to \calI_j (\R^n )$ is given by
$\hat{D}(t,x,\mathbf{v}):=D(\widetilde{c}_k (\mathbf{v})(t))(x)$ which is equal to $D(\widetilde{c}_k^0 (\mathbf{v})(t))(x)$ by the above definition of $\widetilde{c}_k^0$.
Hence the above integral is the same as the value on $\widetilde{c}_k^0$.
\end{proof}

\begin{proof}[Proof of Theorem~\ref{thm:nontrivial} (2)]
By Lemmas~\ref{wheel-like_contribution}, \ref{lem:I(G)=Aut} and \ref{lem:ck0} (1), (3), we have that
\[
	 \pair{\hat{z}_k}{\widetilde{c}_k}-\pair{\hat{z}_k}{\widetilde{c}_k^0}
	= \pair{z_k}{c_k}-\pair{z_k}{{c}_k^0}
	= \pm [\Gamma^{(k)}].
\]
Moreover $\pair{\hat{z}_k}{\widetilde{c}_k^0}=0$ by Lemma~\ref{lem:ck0} (2).
Thus $\pair{\hat{z}_k}{\widetilde{c}_k}$ is equal to $\pm [\Gamma^{(k)}]$, which is not zero by the hypothesis.
This shows that $[\widetilde{c}_k]\in H_{k(n-j-2)}(\femb{n}{j})$ is not zero.

Next we show that $r_*[\widetilde{c}_k]\in H_{k(n-j-2)}(\emb{n}{j})$ is nontrivial when $n,j$ are even and $n\ge 2j$. Consider the following commutative diagram associated with the fibration sequence $ \Omega\Imm(\R^j,\R^n)\stackrel{i}{\to} \femb{n}{j}\stackrel{r}{\to} \emb{n}{j}$:
\[ \xymatrix{
	\pi_{k(n-j-2)}(\Omega\Imm(\R^j,\R^n)) \ar[r]^{i_*} & 
		\pi_{k(n-j-2)}(\femb{n}{j}) \ar[r]^{r_*} \ar[d]^{\overline{H}} &
		\pi_{k(n-j-2)}(\emb{n}{j}) \ar[d]^{H}\\
	& H_{k(n-j-2)}(\femb{n}{j}) \ar[r]^{r_*} &
		H_{k(n-j-2)}(\emb{n}{j}) \\
	}
\]
Here $H$ and $\overline{H}$ are the Hurewicz homomorphisms. The top row is a part of the homotopy exact sequence of the fibration. $H$ and $\overline{H}$ are injective over $\R$ because the component of $\emb{n}{j}$ or $\femb{n}{j}$ ($j\geq 2$) of the standard inclusion is a homotopy associative $H$-space (see \cite[p.263]{MM}). Therefore to show the nontriviality of $r_*[\widetilde{c}_k]$ it is enough to prove the following assertions;
\begin{description}
\item[(a)] $[\widetilde{c}_k]$ lies in the image of $\overline{H}$.
\item[(b)] $r_*\overline{H}^{-1}([\widetilde{c}_k])$ is nontrivial. 
\end{description}
Then (b) and the injectivity of $H$ would imply the result.

Now note that the wheel-like ribbon presentation $P=D\cup B$ in $\R^3$ (Definition~\ref{def:wheel-like-pres}) has the following property:
Let $P'$ be a wheel-like ribbon presentation obtained from $P$ by unclasping the pair $(D_1,B_1)$ as in Figure~\ref{fig:unclasp}.
Then
we can find a 1-parameter family of immersions $\{\varphi_t\}:D^2\to \R^3$, $t\in [0,1]$ such that
(i) $\varphi_0$ is the standard inclusion $\R^j\subset\R^n$,
(ii) $\varphi_t$ restricted to $\partial D^2$ is an embedding for all $t$, and that
(iii) $\varphi_1$ represents $P'$. 
Moreover we may assume that for a base-point $b\in \partial D^2$ and its small neighborhood $U_b$ in $D^2$, it holds that $\varphi_t|_{U_b}=\varphi_0|_{U_b}$ for all $t\in [0,1]$ and thus the connected sum with the standard plane (as in Definition~\ref{def:V_P}) can be done for the entire family.
\begin{figure}
\fig{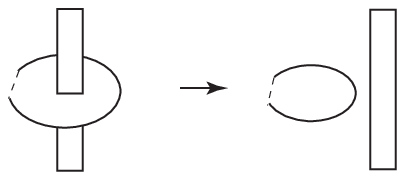}
\caption{}\label{fig:unclasp}
\end{figure}
Then the corresponding family of ribbon $(j+1)$-disks together with embeddings in $\emb{n}{j}$ on its boundaries give a nullhomotopy of a restriction of the map $\widetilde{c}_k:(S^{n-j-2})^{\times k}\to \femb{n}{j}$ to any sub-factor $(S^{n-j-2})^{\times (k-1)}\subset (S^{n-j-2})^{\times k}$. Thus $[\widetilde{c}_k]$ lies in the image of $\overline{H}$ and (a) is proved.

In order to prove (b) we choose a homotopy class $\beta_k\in \pi_{k(n-j-2)}(\femb{n}{j})$ such that $[\widetilde{c}_k]=\overline{H}(\beta_k)$, which exists by (a). $\beta_k$ is nontrivial over $\R$ since $[\widetilde{c}_k]$ is nontrivial over $\R$. Therefore it is enough to prove that in a range $r_*$ on the homotopy group is injective over $\R$.

It is known that $\pi_{l}(\Omega\Imm(\R^j,\R^n))\otimes\R=\pi_l(\Omega^{j+1}V_j(\R^n))\otimes\R$ vanishes for $l\geq 2n-j-6$ (if $n,j$ are even; see \cite[Ch.3, Theorem~3.14]{MT}). Thus, $r_* : \pi_{k(n-j-2)} (\femb{n}{j}) \otimes \R \to \pi_{k(n-j-2)} (\emb{n}{j}) \otimes \R$ is injective if $k(n-j-2) \ge 2n-j-6$. By Proposition \ref{prop:calA} (which will be proved in \S\ref{sec:calA}), $[\widetilde{c}_k]$ can be nontrivial only when $k\ge 3$ (when $n,j$ are even).
It is easy to see that, if $n\ge 2j$, then the above criterion $k(n-j-2)\ge 2n-j-6$ holds for any $k\ge 3$.
\end{proof}

\section{The spaces $\calA_k$}\label{sec:calA}

In this section we discuss the structure of the vector space $\calA_k$.

\subsection{Even codimension case}\label{subsec:n-j=even}

Here we prove the first half of Proposition~\ref{prop:calA}.

\subsubsection{Wheel-type graphs}\label{sssec:wheel-type}
Firstly we introduce the notion of {\em wheel-type graphs} and show that $\calA_k$ is generated by wheel-type graphs
in even codimensional case.

\begin{figure}[htb]
\includegraphics{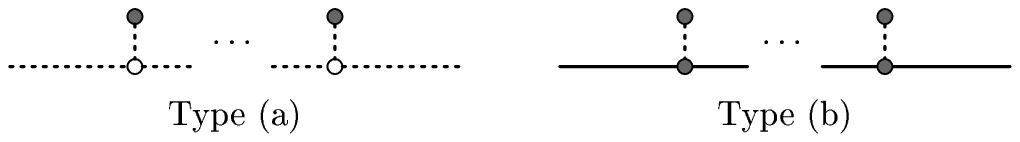}
\caption{Two types of paths}\label{fig:paths}
\end{figure}

\begin{Def}\label{def:wheel}
An admissible 1-loop graph is said to be {\em wheel-type} if it is an alternate cyclic sequence of paths of the form (a) or (b) of Figure \ref{fig:paths}.
A single path may form a loop.
A {\em $k$-wheel} is a wheel-type graph of order $k$ consisting of exactly one path of type (a) (see Figure~\ref{fig:std_lab_wheel}).
We call $\theta$-edges sticking into the paths {\em hairs}.\qed
\end{Def}

\begin{Exa}
Below we show two examples of wheel-type graphs.
\[
 \includegraphics{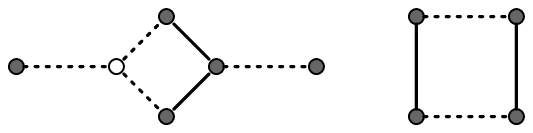}
\]
The left graph consists of one type (a) path and one type (b) path and has two hairs, while the right graph consists of
two type (a) paths and two type (b) paths with no hair.\qed
\end{Exa}

\begin{Lem}\label{lem:all_wheel_type}
In even codimension case, $\calA_k$ is at most one dimensional, possibly generated by the $k$-wheel.
\end{Lem}

\begin{proof}
Let $\Gamma$ be an admissible 1-loop graph, but not wheel-type.
Then $\Gamma$ has at least one tree subgraph $T$ which has $\ge 3$ vertices and shares only one vertex $r$ with the
unique cycle (like the third graph of Example~\ref{ex:1-loop_graph}).
$T$ has one of the following three subgraphs;
\[
 \includegraphics{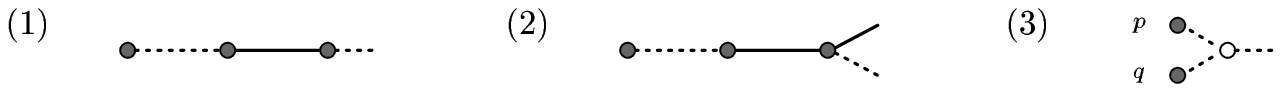}
\]
\underline{Case (1)}.
By the ST relation in Figure \ref{fig:relations_even}, $T$ can be transformed in $\calA_k$ to Case (3).

\noindent
\underline{Case (2)}.
This subgraph is the third one in the ST2 relation (Figure \ref{fig:relations_even}) with the edge $q$ ending at a
univalent vertex.
We can see that the first and the second graphs in the ST2 relation cancel with each other, after the ST and C relations
are applied.
Thus $[\Gamma ]=0\in \calA_k$.

\noindent
\underline{Case (3)}.
Such $\Gamma$ satisfies $[\Gamma ] = -[\Gamma ]$ in $\calA_k$ and hence vanishes, because there is an orientation
reversing automorphism of $\Gamma$ which exchanges $p$ and $q$.

Thus all the graphs which are not wheel-type vanish in $\calA_k$.
As explained in \cite[page 50]{Wa}, by applying relations (Figure \ref{fig:relations_even}), we can transform all the wheel-type graphs to the wheel.
This completes the proof.
\end{proof}

\subsubsection{The case $k\equiv n\equiv j$ modulo $2$}\label{subsubsec:k-n=even}
\begin{figure}[htb]
\includegraphics{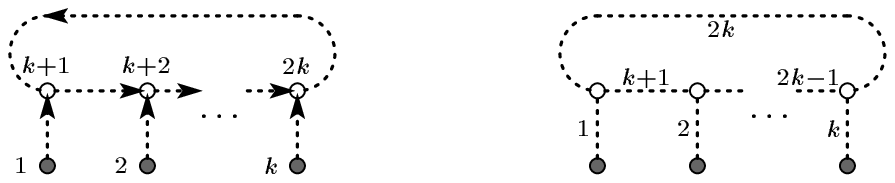}
\caption{Orientations of the $k$-wheel : the cases that $n$, $j$, $k$ odd and that $n$, $j$, $k$ even}
\label{fig:wheel_vanish}
\end{figure}

We can prove that the $k$-wheel vanishes when $k \equiv n \equiv j$ modulo $2$.
Indeed, if we orient the $k$-wheel as in Figure~\ref{fig:wheel_vanish}, then
we can define `reflective' automorphisms $\sigma$ of the $k$-wheel which reverses the orientation as follows:
when $n,j,k$ are odd, $\sigma$ permutes the vertices of the $k$-wheel by
\[
 (1\ \ k)(2\ \ k-1)\dots\Bigl(\frac{k-1}{2}\ \ \frac{k+3}{2}\Bigr) (k+1\ \ 2k)(k+2\ \ 2k-1)\dots\Bigl(\frac{3k-1}{2}\ \ \frac{3k+3}{2}\Bigr)
\]
(whose sign is $(-1)^{k-1}$) and reverses all the $k$ edges on the circle.
When $n,j,k$ are even,
\[
 \sigma :=(1\ \ k)(2\ \ k-1)\dots\Bigl(\frac{k}{2}\ \ \frac{k+2}{2}\Bigr) (k+1\ \ 2k-1)(k+2\ \ 2k-2)\dots\Bigl(\frac{3k-2}{2}\ \ \frac{3k+2}{2}\Bigr)
\]
(whose sign is $(-1)^{k-1}=-1$).
This together with Lemma~\ref{lem:all_wheel_type} proves the following.

\begin{Prop}\label{prop:k=n=j}
If $k\equiv n\equiv j$ modulo $2$, then $\calA_k =\{ 0\}$.\qed
\end{Prop}

\subsubsection{The case $k\not\equiv n\equiv j$ modulo $2$}\label{subsubsec:k-n=odd}
Here we will prove that $\calA_k$ is at least one dimensional if $k\not\equiv n$ modulo $2$.
This will be done by constructing a nontrivial linear map $w_k :\calA_k \to \R$, called a {\em weight system},
for each $k \not\equiv n$ in an analogous way to \cite{Wa}.

\begin{figure}[htb]
\includegraphics{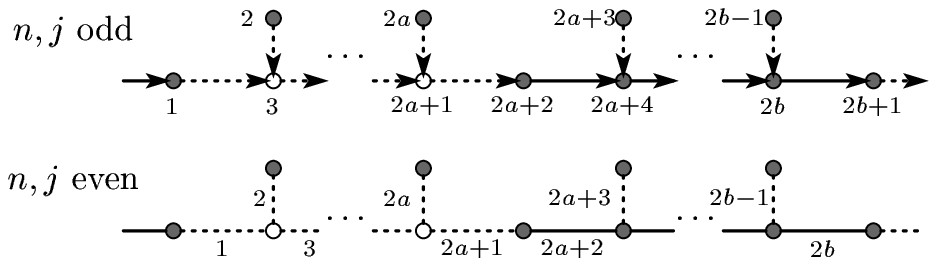}
\caption{Standardly oriented paths (a) and (b)}\label{fig:std_ori_path}
\end{figure}

\begin{figure}[htb]
\includegraphics{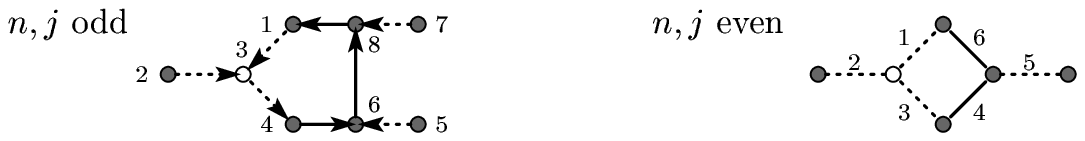}
\caption{Standardly oriented graphs}\label{fig:ex_std_ori}
\end{figure}

\begin{Def}\label{def:std_ori}
A {\em standardly oriented} wheel-type graph is a wheel-type graph oriented as
in Figure~\ref{fig:std_ori_path}.
When both $n$, $j$ are odd, the vertex $1$ is the `first' vertex of a path of type (a), and
when both $n$, $j$ are even, the edge $1$ is the `first' edge of a path of type (a)
(see Figure~\ref{fig:ex_std_ori} for examples).\qed
\end{Def}


\begin{figure}[htb]
\includegraphics{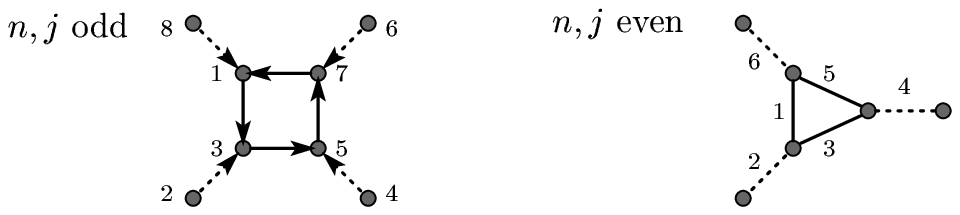}
\caption{Graphs consisting of only one path (b)}\label{fig:graph_one_path(b)}
\end{figure}

\begin{Rem}
There is a unique graph of order $k$ consisting of only one type (b) path.
The standard orientation of the graph is given as in Figure \ref{fig:graph_one_path(b)}.
It is easily checked that this orientation is independent of choices of i-vertex (resp.\ $\eta$-edge)
numbered by $1$.\qed
\end{Rem}

There are
some ambiguities in the definition of the standard orientation; 
the order of the labelling of vertices and edges may be either counterclockwise 
Moreover the definition of a standard orientation depends on the choice of i-vertex/$\theta$-edge numbered by $1$.
But as the name suggests, the standard orientation is uniquely determined.
The proof of the following Lemma is an elementary sign argument.

\begin{Lem}\label{lem:std_ori_welldefined}
Suppose $k \not\equiv n \equiv j$ modulo $2$.
Then any two standard orientations for a wheel-type graph $\Gamma$ of order $k$ are equivalent to each other.\qed
\end{Lem}

For any oriented wheel-type graph $(\Gamma ,\ori (\Gamma ))$ of degree $k \not\equiv n$ modulo 2, define
\[
 w_k (\Gamma ,\ori (\Gamma)):= \varepsilon (-1)^{\sharp \{ \text{hairs of }\Gamma \}}
\]
where $\varepsilon =\pm 1$ is such that $\varepsilon \cdot \ori (\Gamma )$ is equivalent to the standard orientation.
We extend it to a linear map $w_k : \calG_k \to \R$.

\begin{Lem}
When $k \not\equiv n \equiv j$, the map $w_k$ descends to $w_k : \calA_k \to \R$.
\end{Lem}

\begin{proof}
We show that $w_k$ is compatible with the ST relation (Figure \ref{fig:relations_even}) when both $n$ and $j$ are odd.
This relation is represented by the sum of two graphs, which we call $\Gamma_1$ and $\Gamma_2$ respectively
(oriented as in Figure~\ref{fig:relations_even}).
If $\Gamma_1$ is standardly oriented, then so is $\Gamma_2$.
But the numbers of the hairs of $\Gamma_1$ is greater than that of $\Gamma_2$ by one.
Thus we have $w_k (\Gamma_2 ) = -w_k (\Gamma_1 )$ and hence $w_k$ is compatible with the ST relation.
In similar ways we can see that $w_k$ is compatible with all the relations in Figure~\ref{fig:relations_even}.
For the ST2 relation, we may assume the endpoint of the edge labelled by $q$ is univalent since all the graphs here are wheel-type, and then the third graph is zero since it is not wheel-type (see Lemma~\ref{lem:all_wheel_type}).
\end{proof}

\begin{proof}[Proof of Proposition~\ref{prop:calA}, even codimension case]
The case $k\equiv n\equiv j$ modulo $2$ was proved in Proposition~\ref{prop:k=n=j}.
When $k \not\equiv n \equiv j$, we see that $\dim \calA_k \ge 1$, since $w_k (k\text{-wheel})=\pm 1$.
Thus by Lemma~\ref{lem:all_wheel_type}, we have $\calA_k \cong \R$ if $k \not\equiv n \equiv j$.
\end{proof}

\subsection{Odd codimension case}\label{subsec:n-j=odd}
At present we have not determined the structure of $\calA_k$ in odd codimension cases.
Partial descriptions of $\calA_k$ will be given in Propositions~\ref{prop:wheeltype_chord_diagram}, \ref{prop:wheel=0}.
The latter half of Proposition~\ref{prop:calA} will be also proved in Proposition~\ref{prop:A3_odd_codim}.

We call a graph a {\em chord diagram} if it has no e-vertices.
By the defining relations (Figures \ref{fig:relations_odd}, \ref{fig:relation_l}), we can represent every graph
as a sum of chord diagrams in $\calA_k$.
Here we show the following assertion.

\begin{Prop}\label{prop:wheeltype_chord_diagram}
In odd codimension cases, $\calA_k$ is generated by wheel-type chord diagrams.
\end{Prop}

This follows from Proposition~\ref{prop:hair}.
To prove this, we will show the vanishing of chord diagrams with large tree subgraphs introduced in the next two definitions.

\begin{Def}\label{def:feather}
Let $l$ be a positive integer.
A {\em feather} of length $l$ (resp.\ $l +\frac{1}{2}$) is the following subgraph;
\[
 \includegraphics{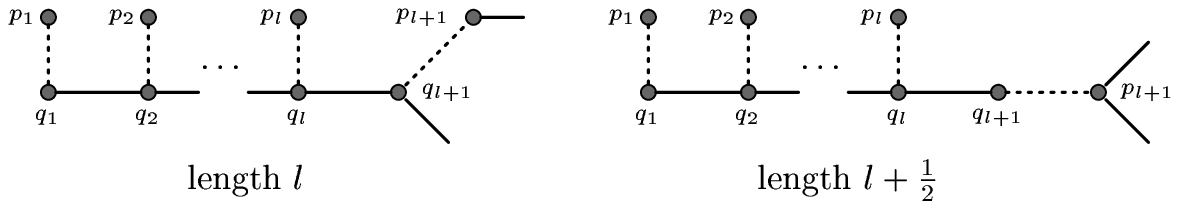}
\]
where $p_1 , \dots ,p_l$ are univalent and $p_{l+1}$ is at least bivalent.
We call the vertex $p_1$ the {\em endpoint} of the feather.\qed
\end{Def}

\begin{Def}\label{def:line}
A {\em straight line} of length $l$ ($l \in \Z_{>0}$) is the following subgraph;
\[
 \includegraphics{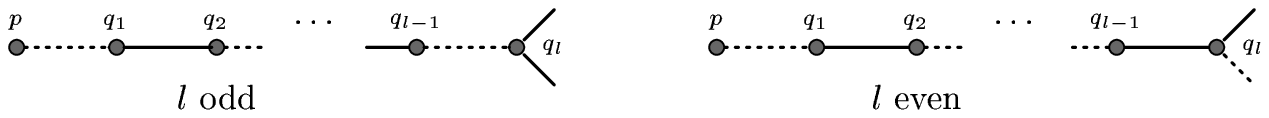}
\]
The vertex $p$ is univalent, $q_1 ,\dots ,q_{l-1}$ are bivalent and $q_l$ is trivalent.
We call the vertex $p$ the {\em endpoint} of the straight line.\qed
\end{Def}

Notice that the straight lines of length $1$, $2$ and $3$ are equal to feathers of length $1/2$, $1$ and $1+(1/2)$,
respectively.
Every univalent vertex is an endpoint of a feather or a straight line.
For example, the vertices $p_2 ,\dots ,p_l$ in a feather are endpoints of straight lines of length $1$.

Below we will prove the following.

\begin{Prop}\label{prop:hair}
In odd codimesion case, any graph can be represented in $\calA_k$ as a sum of chord diagrams all of whose univalent vertices are endpoints of straight lines of length $1$.
\end{Prop}

Any non wheel-type graph must have a subgraph (1) or (2) appearing in the proof of Lemma~\ref{lem:all_wheel_type}, and hence have a straight line of length $>1$.
Hence Proposition~\ref{prop:hair} says that $\calA_k$ is generated by wheel-type chord diagrams,
and completes the proof of Proposition~\ref{prop:wheeltype_chord_diagram}.

The following Lemmas \ref{lem:long_feather_vanish}, \ref{lem:long_line_vanish} and \ref{lem:feather=line} are needed
to prove Proposition \ref{prop:hair}.

\begin{Lem}\label{lem:long_feather_vanish}
If $\Gamma$ has a feather of length $\ge 2+(1/2)$, then $\Gamma =0$ in $\calA_k$.
\end{Lem}

\begin{proof}
The proof for the length $\ge 3$ is as follows;
\[
 \includegraphics{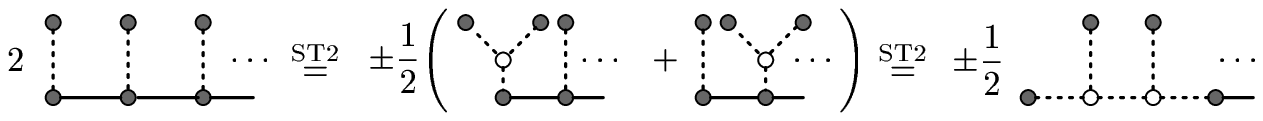}
\]
and the last graph is zero by the IHX relation (see the proof of Lemma~\ref{lem:quasianomalous}, $\Gamma_A$ tree case).

The feather of length $2+(1/2)$ vanishes as follows, again by IHX relation.
\[
 \includegraphics{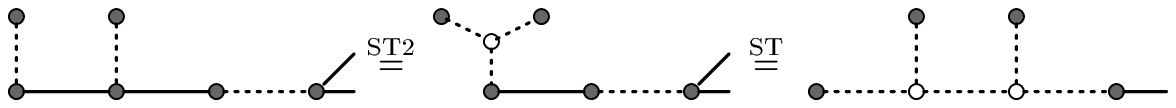}\qedhere
\]
\end{proof}

\begin{Lem}\label{lem:long_line_vanish}
If $\Gamma$ has a straight line of length $\ge 5$, then $\Gamma =0$ in $\calA_k$.
\end{Lem}

\begin{proof}
If the length is at least five, then the straight line contains at least two $\eta$-edges $q_1 q_2$ and $q_3 q_4$
whose endpoints are both bivalent.
Apply the ST relation to $q_1 q_2$ and $q_3 q_4$, then we can transform the straight line to the last subgraph
in the proof of Lemma \ref{lem:long_feather_vanish}.
\end{proof}

\begin{Lem}\label{lem:feather=line}
A straight line of length $4$ is equivalent to the feather of length $2$.
\end{Lem}

\begin{proof}
Apply the ST relation to the $\eta$-edge $q_1 q_2$, and then use the ST2 relation.
\end{proof}

\begin{proof}[Proof of Proposition~\ref{prop:hair}]
Let $\Gamma$ be a chord diagram.
By the above Lemmas \ref{lem:long_feather_vanish}, \ref{lem:long_line_vanish} and \ref{lem:feather=line} and
the fact that the straight lines of length $\le 3$ and the feathers of length $<2$ are equal, we may assume that
all the univalent vertices of $\Gamma$ are endpoints of straight lines of length $\le 4$.

Suppose $\Gamma$ has a straight line of length $>1$.
The straight line of length $4$ can be written by using that of length $3$ as follows;
\[
 \includegraphics{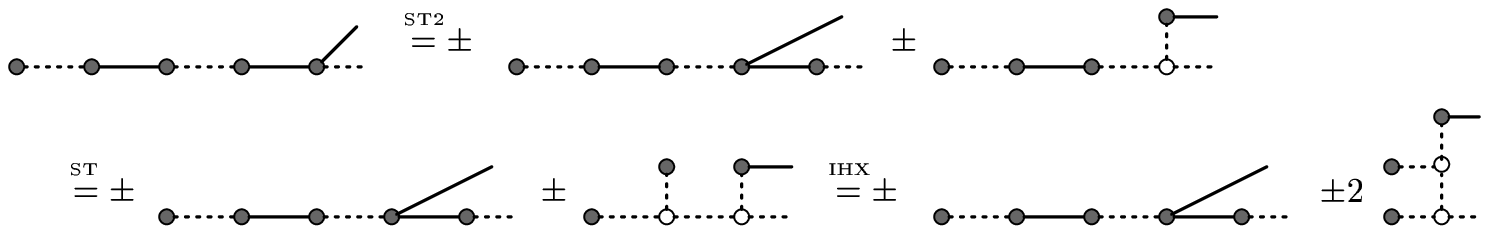}
\]
The last subgraph is equal to that with no univalent vertices by ST relation.

Next we can transform the straight line of length $3$ to a graph with two lines of length $1$;
\[
 \includegraphics{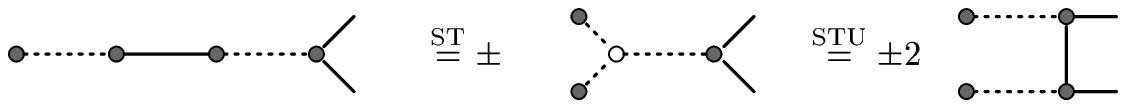}
\]

Lastly the straight line of length $2$ is a sum of a graph with one line of length $1$ and one with no univalent
vertex;
\[
 \includegraphics{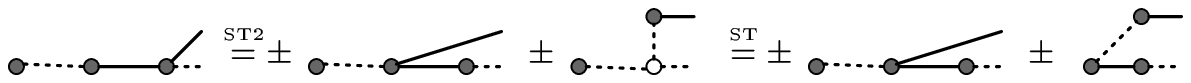}
\]

In such ways as above, we can eliminate all the straight lines of length $>1$.
\end{proof}

We have not yet used the Y relation (Figure \ref{fig:relations_odd}).
The following is a consequence of the ST, STU and Y relations;
\[
 \includegraphics{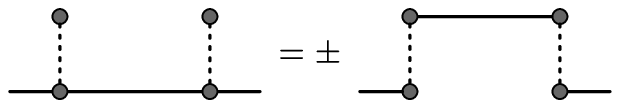}
\]
Thus we can improve Proposition~\ref{prop:wheeltype_chord_diagram} as follows.

\begin{Prop}\label{prop:no_adjacent}
$\calA_k$ is spanned by wheel-type chord diagrams which has no pair of `adjacent' hairs.
\end{Prop}

As a corollary of Proposition~\ref{prop:no_adjacent}, we obtain a very rough, but immediate upper bound of $\dim \calA_k$.
There is exactly one chord diagram with no hair (Figure~\ref{graph_with_no_hair}).
Let $\Gamma$ be a wheel-type chord diagrams with $m>0$ hairs, any two of which are not adjacent to each other.
Then there are $2(k-m)$ bivalent vertices on the cycle of $\Gamma$.
A configuration of hairs determines a partition $2(k-m)=n_1 +\dots +n_m$ (up to cyclic permutations) with all $n_i$'s positive even integers (because there must be even number of bivalent vertices between two non-adjacent trivalent vertices on the cycle).
Then $\dim\calA_k$ is bounded by the number of such partitions.
\begin{Cor}\label{cor:Young_estimation}
We write the number of Young diagrams with $x$ boxes and $y$ rows as $N(x,y)$ (notice that $N(x,y)=0$ if $x<y$).
Then
\[
 \dim \calA_k \le 1+\sum_{1\le m\le \lfloor k/2\rfloor}(m-1)!\cdot N(k-m,m).
\]
For example, we have $\dim\calA_3\le 2$, $\dim\calA_4\le 3$, and so on.
\end{Cor}

The chord diagrams can be obtained by expanding the wheel
by the defining relations.
In this sense the $k$-wheel can be seen as a `source' of the space $\calA_k$.
Thus the next Proposition~\ref{prop:wheel=0} suggests that $\calA_k$ might be rather small in some cases.

\begin{Prop}\label{prop:wheel=0}
\begin{enumerate}
\item
The $k$-wheel vanishes in $\calA_k$ if
(i)  $n$ is even, $j$ is odd, and $k\not\equiv 1$ modulo $4$, or if
(ii) $n$ is odd, $j$ is even, and $k\not\equiv 3$ modulo $4$.
\item
The wheel-type chord diagram which consists of only type (b) paths vanishes if
(i)  $n$ is odd, $j$ is even and $k\not\equiv 1$ modulo $4$, or if
(ii) $n$ is even, $j$ is odd and $k\not\equiv 3$ modulo $4$.
\end{enumerate}
\end{Prop}

\begin{proof}
We prove only (1).
(2) can be proved in a similar way.

Consider the case $n$ is even and $j$ is odd.
Orient the $k$-wheel graph as in Figure \ref{fig:std_lab_wheel} with $\ell$ replaced by $k$;
$(1),\dots ,(3k)$ are $S$-labels, while $1,\dots ,k$ are $T$-labels.
When $k \equiv 3$ modulo $4$, the proof is the same as the argument in \S\ref{subsubsec:k-n=even};
applying the `reflective' permutation which appeared in \S\ref{subsubsec:k-n=even} (whose sign is $-1$) 
to each set $\{ (1),\dots ,(k)\}$, $\{ (k+1),\dots ,(2k)\}$ and $\{ (2k+1),\dots ,(3k)\}$ of the $S$-labels,
we find an orientation reversing automorphism of the $k$-wheel.
Thus the $k$-wheel vanishes.

The proof for even $k$ can be done by applying
the cyclic permutation of $k$ letters (whose sign is $-1$) to each set
$\{ (1),\dots ,(k)\}$, $\{ (k+1),\dots ,(2k)\}$ and $\{ (2k+1),\dots ,(3k)\}$ of the $S$-labels.
The proof for the case $n$ is odd and $j$ is even is similar.
\end{proof}

\begin{figure}[htb]
\includegraphics{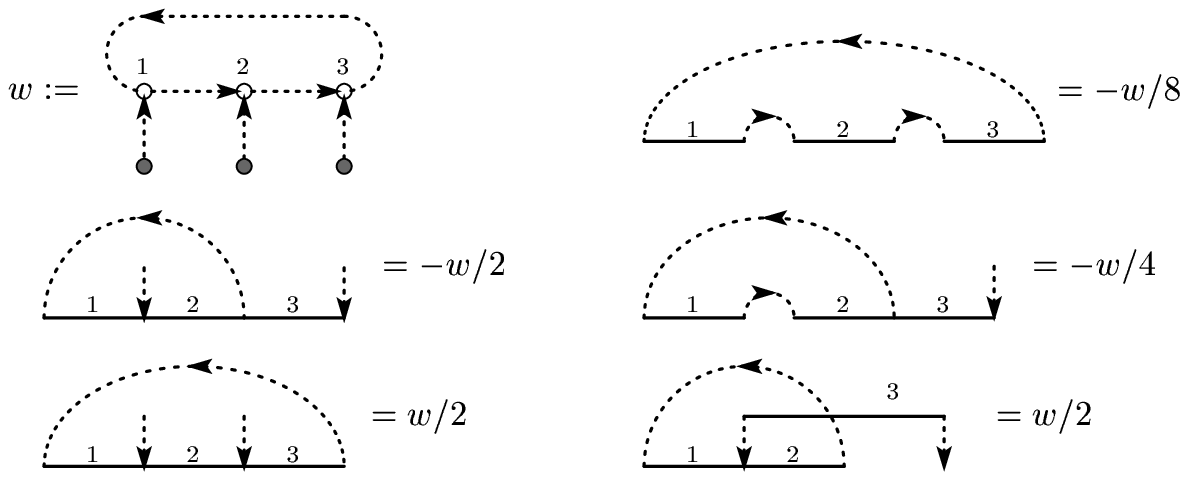}
\caption{All the chord diagrams in $\calA_3$ (i-vertices are omitted)}\label{fig:A_3}
\end{figure}

At present it is difficult to give a lower bound of $\dim \calA_k$, but not impossible if $k$ is small.
Indeed, Figure \ref{fig:A_3} shows all the non-zero chord diagrams in $\calA_3$ which arise from the expansion of the $3$-wheel by the IHX and the STU relations ($n$ odd, $j$ even case).
By solving the system of all possible linear relations among graphs, we can see that all these graphs are equal to the
wheel multiplied by some non-zero constants, and there is no non-trivial relation among these graphs.
Thus we have the following observation.

\begin{Prop}\label{prop:A3_odd_codim}
When $n$ is odd and $j$ is even, the space $\calA_3$ is one dimensional.\qed
\end{Prop}

Since the hexagonal graph (the second graph in Figure \ref{fig:A_3}) does not vanish in $\calA_3$, we obtain a new cohomology class of $\emb{n}{j}$ in odd codimension cases; see Remark~\ref{rem:new_class}.

\providecommand{\bysame}{\leavevmode\hbox to3em{\hrulefill}\thinspace}
\providecommand{\MR}{\relax\ifhmode\unskip\space\fi MR }
\providecommand{\MRhref}[2]{%
  \href{http://www.ams.org/mathscinet-getitem?mr=#1}{#2}
}
\providecommand{\href}[2]{#2}

\end{document}